\documentclass{article}

\usepackage{amsmath}
\usepackage{amssymb}
\usepackage{amsthm}
\usepackage{amsfonts}
\usepackage{amscd}
\usepackage{mathrsfs}
\theoremstyle{plain}

\newtheorem{theorem}{Theorem}[section]
\newtheorem{proposition}[theorem]{Proposition}
\newtheorem{corollary}[theorem]{Corollary}

\theoremstyle{definition}

\theoremstyle{remark}

\newtheorem{remark}[theorem]{Remark}

\numberwithin{equation}{section}

\def\R{ \mathbb R} 
\def\C{ \mathbb C} 
\def\Z{ \mathbb Z} 
\def\T{ \mathbb T} 
\def\zmod#1{ {\mathbb Z}
/#1{\mathbb Z}}
\def\Met{\mathcal Met} 
 
\def\d{ \partial }                  
\def\rttau{ \tau^{1/2} }        
\def\gotimes{\hat{\otimes}} 

\def\Map{\text{\rm Map}}

\def\Aut{\text{\rm Aut}}
\def\Hom{\text{\rm Hom}}
\def\Tr{\text{\rm Tr}}
\def\Det{\text{\rm Det}}
\def\Pfaff{\text{\rm Pfaff}}

\def\sTr{ \text{\rm Tr}_{\text{\rm s}}}

\def\Lie#1{ {\mathfrak #1} }

\def\fpairing#1#2{\langle #1 \wedge #2 \rangle } 
\def\pairing#1#2{\langle #1 , #2 \rangle }

\def\Line#1#2{ {\mathcal L}^{#1}#2 }   
\def\Mod#1#2{ {\mathcal M}_{#1}(#2) } 
\def\Cat#1#2{ {\mathcal C}_{#1}(#2) }  
\def\Mor#1#2{ {\mathcal G}_{#1}(#2) }  

\begin{document} 
\title{Classical Chern-Simons on manifolds with spin structure}
\author{Jerome A. Jenquin}
\date{April 25, 2005}
\maketitle
\begin{abstract}
We construct a 2+1 dimensional classical gauge theory
on manifolds with spin structure whose action is a refinement
of the Atiyah-Patodi-Singer $\eta$-invariant for twisted 
Dirac operators.  We investigate the properties of the Lagrangian
field theory for closed, spun 3-manifolds and compact, spun 
3-manifolds with boundary where the action is interpreted as
a unitary element of a Pfaffian line of twisted Dirac operators.
We then investigate the properties of the Hamiltonian field 
theory over 3-manifolds of the form $\R \times Y$, where
$Y$ is a closed, spun 2-manifold.  From the action we derive 
a unitary line bundle with connection over the moduli stack
of flat gauge fields on $Y$. 
\end{abstract}

By now, Chern-Simons on oriented 3-manifolds is well studied
from various points of view.  On the one hand, Witten 
shows that the quantum theory provides an example
of a Topological Quantum Field Theory, or TQFT \cite{Wi}.  Briefly,
this means that quantum Chern-Simons associates a complex 
number to a closed, oriented 3-manifold; a complex vector space to
a closed, oriented 2-manifold; and algorithms for how to decompose
these associations when cutting along codimension one submanifolds.
(See \cite{A} for a detailed exposition).  His 
arguments, however, pivot around the Feynman Integral, which
is not mathematically well-defined.  Nevertheless, based on these
physical arguments, he is able to show that the 3-manifold invariants 
of this physical TQFT correspond to certain knot invariants.

On the other hand, mathematicians (especially the authors of \cite{BHMV})
construct a TQFT out of Froebenius algebras that they derive from the 
same knot invariants.  Although there is no map between the physical and 
mathematical TQFTs, there is substantial evidence that the two are 
``isomorphic''.  This 
evidence appears at all dimensions of the theory.  In certain cases,
the two TQFTs have been shown to generate the same 3-manifold invariants 
\cite{FG}.  The two TQFTs associate vector spaces of equal dimension to 
a given 2-manifold.  And finally, the two TQFTs associate isomorphic
algebras to the circle.

We see then that Chern-Simons provides a correspondence between a quantized 
field theory and a purely algebro-categorical construction.
Historically, Atiyah deduced the axioms of the latter by abstracting certain 
properties of the Feynman Integral.  But soon after, the algebro-categorical 
TQFTs outpaced the field theoretic TQFTs.  Indeed, mathematicians soon 
discovered knot invariants corresponding to spun 3-manifolds \cite{Li1}, 
\cite{B}.  Out of these Blanchet and Mausbaum construct a new TQFT by 
applying the same categorical techniques developed in \cite{BHMV}.  In this 
new TQFT one finds that the invariants, vector spaces, and algebras 
all depend on the spin structures of the 3-manifolds, 2-manifolds,
and 1-manifolds, respectively.  Blanchet and Mausbaum dub this 
new construction a ``spin-TQFT''.

Chern-Simons provides a correspondence between the the original
``unspun'' TQFTs and a quantized physical theory.  In light of
the construction of Blanchet and Mausbaum one is naturally lead to ask, 
``Where is the physical correspondence to the spin-TQFT?''  
In this paper we put forth an answer of sorts.  

Chern-Simons is a field theory.  In fact, its a gauge theory on 
oriented 3-manifolds.  Only after quantizing do physicists
obtain a TQFT.  What we require is a new field theory on 
3-manifolds; one that incorporates spin structure in a (preferrably)
natural way.  In this paper we construct just such a field theory 
and study its classical aspects.  We dub this new theory 
``spin-Chern-Simons''.  The careful analysis of 
classical spin-Chern-Simons in this paper will lead to an equally 
careful treatment of the quantum theory in the paper \cite{J2}.
Thus, in this paper we propose the ``answer'' to the question above;
while in \cite{J2} we provide evidence that this does
indeed provide the missing correspondence.   

While the ``missing physical correspondence'' question provides a lot of 
motivation for studying spin-Chern-Simons and especially quantum 
spin-Chern-Simons, the classical theory is interesting in its own right.  
Its a field theory that incorporates aspects of classical 
Chern-Simons with geometric index theory, a rather novel feature.  
The action is a spectral invariant associated to operators of Dirac-type; 
and the prequantum line bundle and 
its geometry is Quillen's determinant line with the connection constructed by 
Bismut and Freed \cite{BF1}.  It is the author's hope that this 
approach will one day lead to a solid physical argument for the work of Freed, 
Hopkins, and Teleman \cite{FHT}.  This would require an understanding of  
classical spin-Chern-Simons over 1-manifolds: at present a work in progress. 
For now, we content ourselves with a study of the Lagrangian theory over 
spun 3-manifolds and the Hamiltonian theory over spun 2-manifolds.  Let us
summarize our results.

We construct a classical gauge theory on a 
spin manifold $M$, where $M$ is a compact 3-manifold with spin structure 
$(X, \Sigma )$ for the Lagrangian theory or a compact 2-manifold with 
spin structure $(Y,\sigma )$ for the Hamiltonian theory.  The sigmas label 
the spin structure.  As with most gauge theories, we require some initial data: 
a compact Lie group $G$ and a real, rank zero virtual representation, $\rho$, 
that we call (for now) the ``level''.  

The classical fields are the category of $G$-bundles with connection
over $M$.  That is, the objects are pairs $(P, A)$ where $P \rightarrow M$
is a principal $G$-bundle and $A$ is a connection on $P$; and the
morphisms -- the classical symmetries -- are $G$-bundle morphisms which
cover the identity map of $M$.  So far, this description also applies
to standard Chern-Simons.  We now come to the aspects of spin-Chern-Simons 
that differentiate it from the standard theory.

By choosing a metric on $M$ one naturally has a Dirac operator.
Given any real representation of $G$ one creates a Dirac-{\it like}
operator by coupling a $G$-connection to the Dirac operator
via the representation.  On a 3-manifold with spin structure this operator is
quaternionic and self-adjoint while on a 2-manifold with spin structure it
is complex skew-symmetric.

Over $(X, \Sigma)$ we have a Lagrangian field theory.  To construct
the action, we take any two real representations, $\rho_1$ and $\rho_2$,
such that $\rho_2 - \rho_1 = \rho$.  A $G$-connection,
$A$, induces two Dirac-like operators, $D_1(A)$ and $D_2(A)$,
via $\rho_1$ and $\rho_2$ respectively.  Evaluated at $A$, the
classical action is
$$
{\xi \over 2}(D_2(A)) - {\xi \over 2}(D_1(A)) \quad \text{(mod 1)}
$$
where $\xi $ is the spectral invariant that appears
in the Atiyah-Patodi-Singer index theorem \cite{APS1}.

In general $\xi$ (mod 1) is smooth with respect to smooth
parameters.  Here $\xi /2$ (mod 1) is smooth because the
operators are quaternionic.  Also, note that
the action depends only on the difference $\rho_2 - \rho_1 = \rho$.
The APS index theorem tells us that because $\rho$ is
rank zero, the action is independent of the metric and
that the critical points of the action are exactly the flat
$G$-connections.  Lastly, the action is invariant with respect to
$G$-bundle morphisms, so that it is a well defined function on the
moduli space of $G$-connections.

The action depends on the spin structure $\Sigma$.  Using
the Atiyah-Patodi-Singer {\it flat} index theorem \cite {APS2},
we track this dependence to obtain some useful results.
One result allows us to prove, using cobordism arguments,
that the action only depends on $\rho $ up to the element it
generates in $E^4(BG)$.  Here  $E^*(\cdot )$ is a generalized
cohomolgy theory generated by a spectrum, each element of
which is a twisted product of two Eilenberg-MacLane spaces.
In fact, the degree-four element of the spectrum is the bottom of
the Postnikov tower for $BSO$.  For example, $E^4(BSO_3) 
\cong \Z$ and $E^4(BSU_2) \cong H^4(BSU_2) \cong \Z$ so
that for these groups the spin-Chern-Simons theory has 
integer-valued levels.  

For the Hamiltonian field theory we must consider $G$-connections
over $(Y, \sigma)$.  Given
a real representation $\rho_0$, the $G$-connections over $Y$
provide a family of (complex skew-symmetric) Dirac-like
operators.  Over this family we can take the Pfaffian line bundle
$\Line{\rho_0}{}_Y$ which has a natural unitary structure and
connection \cite{F2}.  In fact, we consider the Pfaffian line of an elliptic 
operator as a {\it graded} line with the grading given by the mod 2 index 
of the operator.  One of the motivations for considering the grading is revealed 
when we discuss the quantum theory in \cite{J2}.

If the level is represented by $\rho = \rho_2 - \rho_1$
then the (graded) line bundle we consider is $\Line{\rho_1}{}_Y \otimes
(\Line{\rho_2}{}_Y)^*$.  $G$-bundle morphisms naturally lift to
a unitary, connection-preserving action on the Pfaffian bundle
so that one gets a unitary line bundle with connection on the
moduli stack of flat $G$-connections $\Mod{G}{Y}$.
This data -- a unitary line bundle with connection over the classical
phase space -- is required when we quantize the Hamiltonian field theory.

If we allow $X$ to have a boundary with $\partial (X,\Sigma )
= (Y, \sigma )$ then by imposing boundary conditions
analogous to those in \cite{DF} the action can still be defined.  In fact,
we consider the exponentiation of the action given by
$$
\tau_X^{1/2}(A) =
{ {\exp{\pi i \xi (D_2(A)) }} \over {\exp{\pi i \xi (D_1(A)) } } }
$$
and we know that $\tau_X^{1/2}(A)$'s
dependence on the boundary conditions naturally identifies it with a
unitary element of the line $\Line{\rho_1}{}_Y\otimes
(\Line{\rho_2}{}_Y)^*|_A$.

If, more generally, $(Y, \sigma ) \subset (X,\Sigma )$, then we 
also have a gluing law that tells us how to factorise $\tau^{1/2}_X(A)$ 
into information on $\Line{\rho_1}{}_Y \otimes (\Line{\rho_2}
{}_Y)^*|_A$.  We point out that the notation ``$\tau_X^{1/2}$'' for the 
exponentiated action is chosen to be consistent with the notation in \cite{DF}.

The outline of this paper is as follows.

In Section 1 we define the action of the theory for $G$-connections
over closed spin 3-manifolds and define the Lagrangian field theory.  
We show that the theory is gauge-invariant and
indepenent of the metric. We compute the Euler-Lagrange equations to see 
that the classical solutions are (gauge group orbits of) flat $G$-connections.  
We also track the action's dependence on the spin structure.  Finally we show 
that the ``levels'' of the theory are actually elements of a certain generalized 
cohomology theory of the the classifying space $BG$.

In Section 2 we define Lagrangian field theory for $G$-connections over 
compact spin 3-manifolds with boundary.  What we find is that, when the 
boundary is non-empty, the action is properly considered to be an element of 
the Pfaffian line for the twisted Dirac operator of the boundary.  To a family 
of $G$-connections the action assigns a section of the Pfaffian line bundle.
With respect to the natural covariant derivative of that bundle we determine
that the section is independent of the metric and derive the same 
Euler-Lagrange equations that we derive is Section 1.  Finally we see that
the action obeys a ``gluing law''.     

In Section 3 we consider the  Hamiltonian field theory over spin 3-manifolds
that are isometric to a spin 2-manifold crossed with an interval.  The 
Euler-Lagrange equations reveal that the space of classical solutions is equivalent to 
the moduli stack of flat $G$-connections over the the 2-manifold.  The action
determines a (Pfaffian) line bundle over the moduli stack and determines the
symplectic structure as well.  We determine that the subset of smooth points of the
moduli stack is the symplectic reduction of the space of $G$-connections with respect 
to the gauge group action.  Finally we point out a functorial 
relationship between the classical Lagrangian and Hamiltonian field theories.

We include an appendix for computations that are necessary but, 
we feel, lie outside the narrative flow of the paper's body.

Finally, the author warmly thanks his thesis advisor, Dan Freed, for introducing
him to this project and for all of his help during its completion.  

\section{The Classical Lagrangian Theory}

\subsection{A review of $\xi$-invariants and the APS-index theorem}

In this section we review some salient points from index theory for
manifolds with boundary as worked out in \cite{APS1, APS2, APS3}.
For the most part -- though much of this discussion applies more
generally -- we focus on 4-manifolds with boundary and closed
3-manifolds.

Let $M$ be an oriented Riemannian 4-manifold with boundary 
$\d M$ and assume that $M$ is spin; that is, assume
$M$ has given a spin structure.  Let $S^{\pm} \rightarrow M$
denote the chiral spinor bundles and consider the chiral Dirac
operator
$$
D_M : C^{\infty}(S^+) \longrightarrow C^{\infty}(S^-)
$$   
between smooth sections of $S^+$ and $S^-$.  If $\d M = \emptyset$ 
then the local form of the Atiyah-Singer Index Theorem tells us that  
\begin{equation}\label{local AS}
\text{index}(D_M) = \int_M 
\left[ \widehat{A} ( \Omega^M ) \right]_{(4)}
\end{equation}
where $\Omega^M$ is the Riemannian curvature two-form
and $\widehat{A}$ is the usual $\widehat{A}$-polynomial so
that
$$
\widehat{A}(\Omega) =
\sqrt{\det  \left( {\Omega / 4 \pi 
\over \sinh \Omega /4 \pi} \right) }.
$$ 
The notation $[ \cdot ]_{(4)}$ is meant to imply the taking the
degree-4 component of the differential form within.

In the case where $\d M$ is not empty, Clifford multiplication
by the (inward pointing) normal vector provides an identification
$S^+|_{\d M} = S^-|_{\d M}$.  We can identify either of these 
with the spinor bundle $S$ of $\d M$.  To obtain something akin 
to the local expression in \eqref{local AS} the authors of 
\cite{APS1} establish global boundary conditions that take 
into account the Dirac operator on the boundary
$$
D_{\d M} : C^{\infty}(S) \longrightarrow C^{\infty}(S).
$$
If we restrict the Dirac operator on M to the subspace of spinor
sections which satisfy these global boundary conditions and the 
geometric data is isomorphic to $[0,1) \times \d M$ near the 
boundary, then the Atiyah-Patodi-Singer (APS) Index Theorem 
offers the expression :
\begin{equation}\label{local APS}
\text{index}(D_M) = \int_M \widehat{A} ( \Omega^M )
- \xi (D_{\d M})
\end{equation}
where $\xi(D_{\d M})$ depends only on the geometric data
of the boundary.  We take a moment now to discuss this 
boundary term as it is the key player in this paper.

The Dirac operator $D$ on any odd-dimensional manifold (such 
as the operator $D_{\d M}$ ) is self-adjoint and elliptic.
Its spectrum is real and discrete so that we may define
$$ 
\eta_{D}(s) = \sum_{\lambda \neq 0}
{ {sign\lambda} \over |\lambda|^s}, \quad Re(s) \gg 0,
$$
where the sum ranges over the nonzero spectrum of $D$.  
One may think of this as a $\zeta$-function regularization 
of the spectral asymmetry
$$
\text{\# of positive eigenvalues of }D - 
\text{\# of negative eigenvalues of }D.
$$ 
It is well established that $\eta_{D}(s)$ converges if $s > 
\dim \d M/2$ and is analytic in $s$; and it has a meromorphic continuation to $s
\in \C $ that is regular at $s=0$.  The value $\eta_{D}(0)$ is the
{\it $\eta$-invariant } of $D$.

What actually appears in the APS-index theorem -- as the reader has already
seen -- is the {\it $\xi$-invariant}:
$$
\xi (D) = { {\eta_{D}(0) + 
\text{dim ker} D } \over 2} 
$$
Under a smooth change of parameters $\xi$ will experience,
at worst, integer jumps so that $\xi$ (mod 1) is smooth.  

All of this also carries over for operators of 
Dirac-type; that is Dirac operators twisted by a unitary connection 
on some vector bundle $E \rightarrow M$.  All that changes 
is the integrand that appears in \eqref{local AS} and 
\eqref{local APS}, which becomes the differential form
$$
\widehat{A}(\Omega^M) ch(\Omega^E).
$$ 
where $\Omega^E$ is the curvature two-form of the connection 
on $E$ and $ch$ is the Chern character polynomial so that
$$
ch (\Omega ) = \Tr \exp (i \Omega / 2 \pi ).
$$
We now take into consideration the dimensions of $M$
and $\d M$.  From the representation theory of $Spin_4$
and $Spin_3$ one has that the bundles $S^{\pm}$ and 
$S$ all have spin-equivariant quaternionic structures
$J^{\pm}$ and $J$ (respectively).  These induce bounded
operators between the $L^2$-sections of the corresponding
spinor bundles, which we also denote by $J^{\pm}$ and $J$.
The Dirac operators are quaternionic by which we
mean
$$
\bar{D}_{\d M} \circ J = J \circ D_{\d M} 
\quad \text{and} \quad
J^- \circ D_M = \bar{D}_M \circ J^+.
$$
There are two upshots to this extra structure: First, the
kernel and cokernel of $D_M$ are each even dimensional
so that index$(D_M)$ is divisible by two.  Second, the 
eigenvalues of $D_{\d M}$ (including zero) occur with 
even multiplicity so that, under smooth changes in the geometric 
parameters, $\xi$ will, at worst experience even valued jumps.  
Thus $\xi /2$ (mod 1) is a smooth function of the geometric 
parameters.  Of course, this discussion more generally applies to 
manifolds of dimensions 4 and 3 (mod 8) due to the Bott 
periodicity of the real Clifford algebras.

To incorporate this refinement with operators of Dirac type we must 
impose some restrictions on the twisting bundle  
$E \rightarrow M$ and its connection.  The discussion above hinged 
upon the existence of quaternionic structures on $S^{\pm}$ and $S$ and 
the compatibility of the respective Dirac operators with these structures.
In order for the bundles $S^{\pm} \otimes_{\C} E$ and 
$S \otimes_{\C} E|{\d M}$ to be quaternionic, $E$ must be real.
In order for the twisted Dirac operators to be compatible, 
the connection on $E$ must be orthogonal.  In these cases 
$\xi /2$ (mod 1) is a smooth function of the geometric and 
twisting parameters.

So far we have only considerd the $\xi$-invariant 
in the context of manifolds with boundary, but this need not always be 
the case.  We may simply look at a closed, compact Riemannian 
3-manifold with spin structure and consider the $\xi$-invariant
of its Dirac operator which is clearly independent of any 4-manifold 
that the 3-manifold might bound.  We do not, however, wish to altogether 
abandon $\xi$'s role as a boundary term in the APS-index theorem.
It is particularly useful when taking the derivative of $\xi$ with
respect to (smoothly) changing geometric data.  To explain this 
we must consider families of Dirac operators.  We follow the 
exposition given in \cite{F1}.

A family of Riemannian 3-manifolds is defined as a fiber 
bundle $\pi : X \rightarrow Z$ with the following structure.
First, the fibers are each diffeomorphic to some given 
compact 3-manifold.  Second, there is a metric structure on 
the relative tangent bundle $T(X/Z)$.  Third, there is  a 
projection $P : TX \rightarrow T(X/Z)$,
the kernel of which is a horizontal distribution on $X$.  We assume 
that $T(X/Z)$ has an orientation and spin structure.  To each 
$z \in Z$ we assign the Dirac operator $D_z$ of the fiber 
$\pi^{-1}(z)$ and this is our family of Dirac operators parametrized 
by $Z$. 

To define a family of twisted Dirac operators parametrized by $Z$
we further require a vector bundle with connection $(E,\nabla) 
\rightarrow X$.  Since $\nabla$ restricts to a connection 
$\nabla^z$ over the fiber $\pi^{-1}(z)$ we can assign to 
$z$ the twisted Dirac operator $D_{E_z}$.  In the context of 
families $\xi$ (mod 1) gives us a smooth function $z \mapsto 
\xi (D_{E_z})$ (mod 1) on $Z$.  

The APS-index theorem can be 
used to track changes in $\xi$ (mod 1) along any smooth path in $Z$. 
In the untwisted case one sees that, for infinitesimal changes, 
$$
d \xi =  \left[ \int_{X/Z} 
\widehat{A}(\Omega^{X/Z}) \right]_{(1)}
$$  
where $\Omega^{X/Z}$ is the curvature two form of the
relative Levi-Civita connection.  The notation is meant to
imply that we are integrating the form $\widehat{A}
(\Omega^{X/Z})$ over the fibers and considering only
the degree-1 component on $Z$.  There are, of course, corresponding 
expressions for the $\xi /2$-invariant and families of twisted 
Dirac operators given by
\begin{equation}\label{difofxi}
\frac{d \xi}{2} =  \frac{1}{2}\left[ \int_{X/Z} 
\widehat{A}(\Omega^{X/Z}) ch (\Omega^E) \right]_{(1)}
\end{equation}   

Before we end this section we would like to point out how
the $\xi$-invariants of twisted Dirac operators behave 
under direct sums of the twisting data.  Let 
$(E_1, \nabla^{E_1})$ and $((E_2, \nabla^{E_2})$
be two vector bundles with unitary connections and let
$D_{E}$ denote the Dirac operator twisted by $(E,\nabla^E)$.
Then, as is easy to see from its definition, 
\begin{equation}\label{addition}
\xi (D_{E_1 \oplus E_2} ) = 
\xi (D_{E_1}) + \xi (D_{E_2})
\end{equation}
so that $\xi $ behaves additively under direct sums.
This motivates an obvious definition for the 
$\xi $-invariant of a Dirac operator twisted by 
{\it virtual} vector bundles with connection.  The
``virtual'' here just means formal differences.  Thus,
for the virtual vector bundle with connection
$(E_1,\nabla^{E_1}) - (E_2,\nabla^{E_2})$,
we define
\begin{equation}\label{subtraction}
\xi (D_{E_1 - E_2}) =
\xi (D_{E_1}) - \xi (D_{E_2})
\end{equation}
Of course, this also applies to the refined $\xi /2$-invariant
and virtual real vector bundles with connection.

\subsection{The data of the classical gauge theory}
In this section we will describe the data that goes into 
defining our classical field theory.  

To begin, fix a closed, compact,  Riemannian 3-manifold $X$
with spin structure.  Let $D$ denote the corresponding Dirac operator.  
It is well known that any oriented 3-manifold admits a spin 
structure so that this imposes no constraints on $X$.  Now we 
define our classical gauge theory over $X$.  

Any gauge theory requires some structural data.  Our gauge theory 
requires a compact Lie group $G$ and a virtual real
representation $\rho = \rho_2 - \rho_1$ of rank zero.  For now,
we call $\rho$ the ``level'' of our theory.  

Given this structural data we take our ``fields'' to be objects
in the category $\Cat{G}{X}$ of principal $G$ bundles 
with connection over $X$.  Throughout we refer to the 
objects in this category as {\it G-connections} on $X$ or just
{\it connections} when $G$ or $X$ is understood.  The ``symmetry'' 
group $\Mor{G}{X}$ of this theory is taken to be $G$-bundle 
morphisms which cover the identity of $X$.  
That is, we consider any two connections to be physically equivalent
when there is a morphism (covering $id_X$) that maps one to 
the other.

Before we define the action of our field theory we say 
a few words about associated vector bundles and establish
some notation.  This will facilitate the definition of our action
and much of what follows.  Given
a pair $(P,A)$ -- where $P \rightarrow X$ is a $G$-bundle
and $A$ is a connection on $P$ -- the associated bundle
construction assigns to any representation $\rho' : G
\rightarrow GL(V)$ a vector bundle with connection 
$$
(\rho' P,\rho' A) =
(P \times_{\rho'} V , d + \dot{\rho'}(A) ).
$$
Here $d$ is the exterior derivative on $P$ and $\dot{\rho'}(A)$ 
is the $\Lie{gl}(V)$-valued one form obtained from $A$
and the Lie algebra homomorphism induced by $\dot{\rho'} : 
\Lie{g} \rightarrow \Lie{gl}(V)$. 
The connection $d + \dot{\rho'}(A)$ acts on equivariant maps $P 
\rightarrow V$ which are, afterall, the sections of 
the associated bundle $\rho' P$.
Likewise, this construction assigns to any virtual 
representation (like the level $\rho$ above) a 
virtual vector bundle with connection, the formal difference 
of two representations becoming the formal difference of the
two corresponding vector bundles with connection.  

That being said, we define our action to be the map
\begin{align*}
\Cat{G}{X} & \longrightarrow \R /\Z  \\
(P,A)    & \longmapsto \xi /2 (D_{\rho A} ) 
\quad \text{(mod 1)}
\end{align*}  
where $D_{\rho A}$ is the Dirac operator on $X$ twisted by
the virtual vector bundle with connection $(\rho P,\rho A)$.

Let $\T$ denote the elements of $\C$ with modulus 1. 
As our action takes values in $\R / \Z$ it is perhaps just 
as natural (if not more natural as we shall see later)
to take the action to be the $\T$-valued exponential of the 
$\xi /2$-invariant.  In \cite{DF} Dai and Freed analyze
the exponentiated $\xi$-invariant 
$\tau_X(D) = \exp{2 \pi i \xi (D)} $
and following their notation we define
$$
\rttau_X (D) = \exp{2 \pi i \xi (D) / 2}
$$
the square root being well define for the reasons
given in the previous section.  From here on out
our action will be the $\T$-valued map
\begin{align*}
\Cat{G}{X} & \longrightarrow \T \\
(P,A)    & \longmapsto \rttau_X (D_{\rho A} ). 
\end{align*}  
Though, if the 3-manifold $X$ is understood or 
peripheral, it is sometimes dropped from the notation.

Given the structure and the action described above
we call this classical field theory ``spin-Chern-Simons'' 
at level $\rho$.  This name is purposely chosen to 
invoke the memory of classical Chern-Simons theory, where the 
action is defined to be the Chern-Simons 3-form associated to a 
G-connection \cite{CS}, \cite{F2}.

\subsection{Functoriality of $\rttau$}
In this section we enumerate some of the features of our action.

The structure of this next proposition emulates Proposition 2.7 of
\cite{F2}.

\begin{proposition}\label{functorial}
Let $X$ be a closed oriented 3-manifold with Riemannian
metric $g$ and spin structure $\sigma $.  Then the
spin-Chern-Simons action
$$
\rttau_X: \Cat{G}{X} \longrightarrow \T
$$  
satisfies the following properties:
\begin{enumerate}
\item (Functoriality) If $\varphi: P' \rightarrow P$ and 
$F: S'_{g', \sigma'} \rightarrow S_{g, \sigma }$ are any 
G-bundle and spinor bundle morphisms (respectively) covering an 
orientation and spin structure preserving isometry $f: (X',g') 
\rightarrow (X,g)$, and $A$ is a connection on P, then
$$
\rttau_{X'} (D'_{\rho (\phi^*A)}) = \rttau_{X} (D_{\rho A})
$$
where $D'$ is the Dirac operator on $X'$.
\item (Orientation) Let $-X$ denote $X$ with the opposite 
orientation. Then 
$$
\rttau_{-X}(D_{\rho A}) = 
\overline{ \rttau_X(D_{\rho A}) }.
$$
\item (``Additivity'') If $X = X_1 \sqcup X_2$ is a disjoint union,
and $A_j$ are connections over $X_j$, then
$$ 
\rttau_{X_1 \sqcup X_2}(A_1 \sqcup A_2) =
\rttau_{X_1}(A_1) \cdot \rttau_{X_2}(A_2).  
$$
\end{enumerate}  
\end{proposition}

It follows from (1) that the action is invariant under the symmetry
group $\Mor{G}{X}$ of the theory.  Or what is equal, there is an 
induced action
$$
\rttau_X :\overline{\Cat{G}X} \longrightarrow \T
$$
where $\overline{\Cat{G}X}$ is the set of equivalence classes.

\begin{proof}
To prove (1) we point out that
$$
D'_{\rho (\varphi^*A)} = 
\Phi_{F,\varphi }^{-1} \circ D_{\rho A}
\circ \Phi_{F,\varphi } 
$$
where $\Phi_{F,\varphi }$ denotes the induced unitary 
map on sections given by
$$
s' \longmapsto (F \otimes \rho \varphi ) \circ s' \circ f^{-1}.
$$
Thus $D'_{\rho \varphi^*A}$ and $D_{\rho A}$ have the
same eigenvalues and eigenvalue degeneracies so that their
respective $\rttau$-invariants are equal.

To prove (2) we point out that on oriented spin 3-manifolds -- 
and oriented odd dimensional spin manifolds in general --
the orthonormal volume form provides a covariantly
constant, spin-equivariant isomorphism between spinor bundles 
$\omega : S_X \rightarrow S_{-X}$.  We denote the induced
unitary map on sections by $\omega$ as well.  A standard argument
shows that $D_{-X} \circ \omega = -\omega \circ D_X$.  Thus,
if $\lambda $ is an eigenvalue of $D_X$, $-\lambda$ will be an
eigenvalue of $D_{-X}$.

The proof for (3) is easy and we leave it to the reader.  
\end{proof}

The $\rttau$-invariant also behaves well under extension of the 
structure group.  We explain what we mean.  An inclusion of 
Lie groups $i: G \hookrightarrow G'$ induces an inclusion of
principal bundles $i_P : P \hookrightarrow P'$, where $P$ is 
any $G$-bundle and $P' = P \times_i G'$.  If $A$ is a connection 
on $P$ then there is a natural extension to a connection $A'$ on 
$P'$ determined by $i (A) = i_P^*(A')$ (cf \cite{KN} ).  This
is the situation we consider in the next propoposition.

\begin{proposition}\label{ext of grp}
Let $\rho$ be a virtual representation of $G'$.
For any $G$ connection $A$ over $X$, its extension $A'$ to a 
$G'$ connection satisfies $\rttau (D_{(\rho \circ i) A}) =
\rttau (D_{\rho A'})$.  
\end{proposition} 

Another way to express Proposition \ref{ext of grp} is 
to say that the extension functor $i_X : (P,A) \mapsto 
(P ', A')$ fits into the following commutative 
diagram:
$$\begin{CD}
\Cat{G}X              @>i_X>>     \Cat{G'}X \\
@V{\rttau}VV                         @VV{\rttau}V \\
\T                @>id_{\T }>>           \T
\end{CD}$$
 
 \begin{proof}
There is a natural isomorphism between the associated virtual vector
bundles $\rho P' \rightarrow (\rho \circ i)P$ which sends the associated
virtual connection $\rho A'$ to $(\rho \circ i) A$.  Thus the twisted 
Dirac operators $D_{\rho A'}$ and $D_{(\rho \circ i) A}$ will have the
same eigenvalues and eigenvalue degeneracies and so the same 
$\rttau$-invariants.   
\end{proof}
 
\subsection{Dependence of $\rttau$ on smooth parameters}

To facilitate the discussion which follows we establish some conventions and
notation.  Any real representation $\rho'$ of a compact lie group $G$ generates an
Ad-invariant bilinear form on the lie algebra ${\frak g}$:
\begin{equation}\label{p1 form}
\pairing{\eta_1}{\eta_2}_{\rho'} = 
-{1 \over 8\pi^2} \Tr ( \rho' (\eta_1 ) \rho' (\eta_2)  )
\end{equation}
\begin{remark}
It is easy to check  that if $[\eta_1, \eta_2] = 0 $ and $ \exp (\eta_1 ) =
\exp( \eta_2 ) = 1_G$ then $\pairing{\eta_1}{\eta_2}_{\rho'}$ is 
$\Z$-valued.  This provides some motivation for the normalization.
Another motivation comes from Chern-Weil theory.  
If $P$ is a principal $G$-bundle and $\rho' P$ is the associated 
vector bundle then the first Pontryagin class $p_1(\rho' P)$ can be 
represented in de Rham cohomology by the 4-form
$$
\fpairing{ \Omega^A}{\Omega^A}_{\rho' }.
$$
for any connection $A$ on $P$.
Notice that if  $\rho = \rho_2 \oplus \rho_1$ then
$$
\pairing{\eta_1}{\eta_2}_{\rho } = \pairing{\eta_1}{\eta_2}_{\rho_2 } +
\pairing{\eta_1}{\eta_2}_{\rho_1 }
$$
For this reason, if we consider a {\it virtual} representation
$\rho = \rho_2 - \rho_1$ we define its Ad-invariant bilinear form 
$\pairing{}{}_{\rho }$ by 
$$
\pairing{\eta_1}{\eta_2}_{\rho } =
\pairing{\eta_1}{\eta_2}_{\rho_2 } - 
\pairing{\eta_1}{\eta_2}_{\rho_1 }
$$
\end{remark}

As it is defined the action seems to depend on the 
Riemannian structure on $X$.  But, in fact we have the 
following:
 
\begin{proposition}\label{d wrt metric}
If $\rho$ is a real rank zero virtual representation then
$\rttau (D_{\rho A})$ is independent of the metric.
\end{proposition}

\begin{proof}
If $\Met (X)$ denotes the (convex) space of metrics 
on $X$ then we take the obvious family 
$\Met (X) \times X \rightarrow \Met (X)$.
We fix the pair $(P,A)$ and compute the differential
over $\Met (X)$ using the twisted version of \eqref{difofxi}:
$$
d \rttau = \rttau \left[ \pi i \int_{X} 
\widehat{A} (\Omega )
ch(\rho \Omega^A ) \right]_{(1)}.
$$
where (for typographical reasons) the script-less 
$\Omega$ is standing in for the Riemannian 
curvature of the relative tangent bundle $T( (\Met (X) 
\times X) / \Met (X) )$.  Unraveling the polynomials
$\widehat{A}$ and $ch$ we get a more explicit
expression for the integrand
$$
\widehat{A}(\Omega )
ch(\rho \Omega^A ) =
\text{rank}(\rho ) \widehat{A}
(\Omega )_{(4)}
+ \fpairing{\Omega^A}{\Omega^A}_{\rho }
+ \text{higher degree terms}.
$$ 
Since $\rho$ is rank zero the first term goes
away.  The second term -- since $A$ is fixed --
is the pullback (by the vertical projection
$\Met (X) \times X \rightarrow X$) of a 
4-form on a 3-manifold and is thus identically zero.  
\end{proof}

Thus the Lagrangian field theory is independent of
the metric on $X$.  We would now like to see how
$\rttau$ behaves with respect to infinitesimal 
changes in the connection (the only remaining smooth
parameter in our theory).  Let $\Cat{}{P}$ denote the
affine space of connections on a principal $G$-bundle
$P \rightarrow X$.  Any smooth variation of connections
will occur within such a space. 

\begin{proposition}\label{d wrt conn}
Let $A_t$ be a path in $\Cat{}{P}$ and let $\Omega_t$
denote the curvature of the connection $A_t$, $\dot{A_t}$
the tangent vector along the path.  Then
$$
 \frac{d}{dt}\rttau ( A_t) = 2 \pi i \cdot dt \cdot
\rttau (A_t) \int_X 
\fpairing{\Omega_t}{\dot{A_t}}.
$$
\end{proposition}

\begin{proof}
The connections $A_t$ define a family of twisted Dirac 
operators parametrized by $[0,1]$.  Or what is the same,
they give us a single connection on the $G$-bundle
$P \times [0,1]$.  The curvature of this single connection
is $dt \wedge \dot{A_t} + \Omega_t$. As before, we use 
the twisted version of \eqref{difofxi} to compute
\begin{align*}
\frac{d}{dt} \rttau (A_t) & = 
\rttau (A_t ) \left[ \pi i \int_X 
\fpairing{(dt \wedge \dot{A_t} + \Omega_t )}
{(dt \wedge \dot{A_t} + \Omega_t )} \right]_{(1)} \\
& = \pi i \cdot dt \cdot \rttau (A_t) \int_X 
2 \fpairing{\dot{A_t}}{\Omega_t}
\end{align*}
proving the proposition. 
\end{proof}

\begin{remark}\label{mod space1}
Assuming $\pairing{}{}_{\rho}$ is nondegenerate, this 
proposition implies that $d \rttau |_A = 0$ if and only if
$\Omega^A = 0$; that is, if and only if $A$ is flat.  This
is the content of the Euler-Lagrange equation
\begin{equation}\label{eul-lag} 
\Omega^A = 0
\end{equation}
which is first order.  Since $\rttau_X$ is invariant with
respect to gauge transformations $\Mor{G}{X}$, so too is
the space of solutions to the Euler-Lagrange equation.
This is also obvious from \eqref{eul-lag}.  We let 
$$
\Mod{G}{X} \subset \overline{\Cat{G}{X}}
$$
denote the space of equivalence classes of solutions
to \eqref{eul-lag}.  We will have more to say about
$\Mod{G}{X}$ when we discuss the Hamiltonian
field theory in Section 1.3. 
\end{remark}

\subsection{Dependence of $\rttau$ on the spin structure}

We now track the dependence of the action on the spin structure
$\sigma$ assigned to $X$.  This piece of the data is discrete.  
In fact, if $spin(X)$ denotes the set of equivalence classes of spin 
structures on $X$, then $spin(X)$ is affine over the vector space 
$H^1(X ; \zmod2)$.  This vector space is in one-to-one 
correspondence with equivalence classes of flat, orthogonal line 
bundles on $X$ so that, throughout the rest of the paper, we often 
identify a flat line bundle with its first Stiefel-Whitney class.

If $\ell$ is a flat line bundle and $S_{\sigma}$ is the spinor
bundle associciated to the spin structure $\sigma$, then recall
that we can identify $S_{\sigma + \ell }$ with $S_{\sigma} 
\otimes \ell$ -- we assume the metric is fixed.  The Dirac
operator associated to $\sigma + \ell$ is therefore the 
Dirac operator associated to $\sigma$ twisted by $\ell$.  

To track the dependence of the action on $\sigma $ we fix
the metric and the pair $(P,A)$.  From the dicussion in the 
previous paragraph we would like to compute the ratio
\begin{equation}\label{quad form}
q_{\sigma }(\rho P, \ell ) = 
\frac{ \rttau ( D_{\ell \otimes \rho A} ) }
{\rttau( D_{\rho A } )} =
\rttau ( D_{(\ell -1) \otimes \rho A} )
\end{equation}  
effectively comparing the $\rttau$-invariants for the 
twisted Dirac operators on $S_{\sigma } \otimes \rho P$
and $S_{\sigma } \otimes \ell \otimes \rho P$.  

We point out that, since $\Omega^{\rho A \otimes \ell}
= \Omega^{\rho A}$, Proposition \ref{d wrt conn} implies
that $q_{\sigma }$ is independent of the connection
$A$.  Along with that Proposition \ref{functorial} implies that
$q_{\sigma }$ depends only on the topological type
of $P$ as $\rttau$ is invariant under bundle morphisms.
Furthermore, \eqref{addition} and \eqref{subtraction} imply 
that $q_{\sigma }$ depends only on the element of 
$KO(X)$ represented by the virtual vector bundle $\rho P$.
Thus $q_{\sigma }$ depends only on discrete 
topological parameters.  Indeed, it has an entirely 
KO-theoretic interpretation given by the Atiyah-Patodi-Singer
(APS) Flat Index Theorem \cite{APS2}.   We state the 
nature of $q_{\sigma }$ in the following proposition.

\begin{proposition}\label{properties of q}
Let $E$ be a rank zero element of $KO(X)$. 
\begin{enumerate}
\item If $w_1(E) = 0$ and $w_2(E) = 0$ then 
$q_{\sigma }(E, \cdot ) \equiv 1$.
\item If $w_1(E) = 0$ then $q_{\sigma }(E, \ell ) =
(-1)^{w_2(E) \smile \ell }$.
\item In general, $q_{\sigma }(E, \cdot )$ is a 
$\zmod4$-valued quadratic
refinement of the bilinear form 
$$
B_{E}( \ell_1, \ell_2) =
{q_{\sigma }(E,\ell_1 \otimes \ell_2) \over 
q_{\sigma }(E,\ell_1)q_{\sigma }(E,\ell_2)}
= (-1)^{w_1(E) \smile  \ell_1 \smile \ell_2}.
$$
\end{enumerate} 
\end{proposition} 

\begin{proof}
To prove (1) we point out that $q_{\sigma }$ depends
linearly on $E$.  Indeed, if $w_1(E) = 0$ and $w_2(E) = 0$
Propostion \ref{ko of x} in the appendix tells us 
that $E = 0$ so that $q_{\sigma }(E, \cdot ) \equiv 1$.

To prove (2) we first point out that, if $w_1(E)=0$ then an
easy computation shows that $w_1(2E) = 0$ and $w_2(2E) = 0$
so that $2E = 0$.  Thus $q_{\sigma }(E,\cdot )$ takes
values in $\zmod2^{\times }$.  Recall that, according the Flat Index 
Theorem, $q_{\sigma }$ is a cobordism invariant. (One 
can easily show this using the twisted version of 
\eqref{local APS}).  Thus we have a homomorphism
$$
q_{\sigma } : \Omega_3^{spin} (BSO \times B(\zmod2))
\longrightarrow \zmod2^{\times}
$$   
where $\Omega_3^{spin} (BSO \times B(\zmod2))$ is the
cobordism group of smooth 3-manifolds with spin structure,
oriented rank zero element of $KO$, and flat line bundle.
Even better, assuming the validity of (1), it is clear that this
homomorphism factors through $\Omega_3^{spin} (BSO 
\wedge B\zmod2)$ which, 
in Proposition \ref{cobordism grps} of the appendix,
we show is isomorphic to $\zmod2$.  But in Proposition 
\ref{q on 3torus} we show that, on the 3-torus, we have exactly the 
formula
$$
q_{\sigma }(E, \ell ) =
(-1)^{w_2(E) \smile \ell }
$$
called for in (2) of this proposition.  Since Stiefel-Whitney
numbers are also cobordism invariants \cite{Hu} and
the structures on the 3-torus clearly generate 
$\Omega_3^{spin} (BSO \wedge B\zmod2)$ this gives us
(2) for general compact 3-manifolds.

To prove (3) we decompose $E = (\ell - 1) + F$ where 
$w_1(E) = \ell$ and $F \in \widetilde{KO}(X)$ is 
such that $w_1(F) = 0$ and $w_2(F) = w_2(E)$.  Then 
$$
B_{E}( \ell_1, \ell_2) = B_{\ell -1}(\ell_1,\ell_2)
\cdot B_{F}(\ell_1,\ell_2) = B_{\ell -1}(\ell_1,\ell_2)
$$
since $B$ is linear with respect to $E$ and is zero if 
$w_1(E) = 0$.  Now, rearranging some the the factors,
we see that 
\begin{align*}
B_{\ell -1}(\ell_1,\ell_2) & =
\rttau (D_{(\ell -1) \otimes 
(\ell_1 \otimes \ell_2 \oplus 1 - \ell_1 - \ell_2)}) \\
& = q_{\sigma }(\ell_1 \otimes \ell_2 
\oplus 1 - \ell_1 - \ell_2, \ell ) \\
& = (-1)^{\ell_1 \smile \ell_2 \smile \ell}
\end{align*}
where the last equality follows from (2) and the fact 
that $w_1(\ell_1 \otimes \ell_2 \oplus 1 - \ell_1 - 
\ell_2) = 0$ and $w_2(\ell_1 \otimes \ell_2 
\oplus 1 - \ell_1 - \ell_2) = \ell_1 \smile \ell_2$.
\end{proof}

\subsection{Redefining the ``level''} 

As it is defined above the level is a real, rank zero, virtual 
representation $\rho = \rho_1 - \rho_2$.  Even at first 
glance one sees that we need not be so specific.  In particular, 
fixing all other parameters -- spin structure, connection, 
metric -- $\rttau (D_{\rho A})$ only depends on the 
equivalence class represented by $\rho$ in the rank zero 
representation ring $\widetilde{RO}(G)$.  If we look 
a little harder, in fact, we can see that even this is too specific 
and that the level need not be so refined.  We elaborate on this next.
 
In Section A.3 of the appendix we define a generalized 
cohomology theory $X \mapsto E^{\bullet }(X)$.
As discussed in the appendix, a virtual representation
$\rho$ determines a class $\lambda (\rho ) \in E^4(BG)$.
Our claim is that $\rttau (D_{\rho A})$ really only
depends on the class $\lambda (\rho )$.  Or, what's equivalent,
that the map 
\begin{align*}
\widetilde{RO}(G) & \longrightarrow \T \\
\rho                          & \longmapsto  \rttau (D_{\rho A})
\end{align*}   
factors through the map $\lambda : \widetilde{RO}(G)
\rightarrow E^4(BG)$.  This is exactly what is implied 
by the following proposition.

\begin{proposition}\label{level}
If $\lambda (\rho ) = 0$ then $\rttau (D_{\rho A}) = 1$
for all $G$-connections $A$.
\end{proposition}

\begin{proof}
Let $\rho = \rho_1 - \rho_2$.  Without loss of generality
we can assume that both $\rho_1$ and $\rho_2$ map $G$ into
$SO(V)$.  Recall -- again, from the appendix -- that $E^4(BG)$
fits into an exact sequence
$$
0 \rightarrow H^4(BG) \rightarrow E^4(BG)
\rightarrow H^2(BG;\zmod2)
$$ 
where the third homomorphism is induced by the 2nd Stiefel-Whitney
class $w_2 \in H^2(BSO)$.  Given the hypothesis of the 
proposition this implies that $w_2  (\rho ) 
= 0$ so that, for any $G$-bundle $P$, $w_2 ( \rho P ) = 0$.
Together with Proposition \ref{properties of q}, this implies
that $\rttau (D_{\rho \cdot})$ is independent of the spin 
structure.
  
Yet another fact from Section A.3 is that the first Pontryagin class
$p_1 \in H^4(BSO)$ induces a homomorphsim $E^4(BG) 
\rightarrow H^4(BG)$.  Together with the hypothesis this 
implies that, for any $G$-bundle $P$, $p_1 (\rho P) = 0$.
Furthermore, from the Chern-Weil isomorphsim 
$$
H^{\bullet} (BG) \otimes \R \cong 
Sym^{\bullet }(\Lie{g}^*)^{Ad G},
$$ 
we see that $\pairing{}{}_{\rho } = 0$

This last revelation implies that $\rttau (D_{\rho \cdot})$
is a cobordism invariant.  Indeed, if all of the parameters on $X$
-- the pair $(P_X,A_X)$, metric $g_X$, and spin structure $\sigma_X$
-- bound corresponding parameters on a 4-manifold $M$ -- 
$(P_M,A_M)$, $g_M$, $\sigma_M$ -- then according to 
\ref{local APS}
$$
\rttau (D_{\rho A}) = \exp \pi i \left(
\int_M \pairing{\Omega^{A_M}}
{\Omega^{A_M}}_{\rho } \right ) = 1. 
$$ 
Since there is no obstruction to extending the smooth
parameters, $\rttau (D_{\rho (\cdot )})$ induces a 
homomorphism $\Omega_3^{spin}(BG) \rightarrow \T$,
where $\Omega_3^{spin}(BG)$ is the cobordism group of 
3-manifolds with spin-structure and a $G$-bundle.  In fact,
this homomorphism factors through $\Omega_3^{spin}(BSO)$.  

Now we borrow one last fact from the appendix.  In 
Section A.1 we find that $\Omega_3^{spin}(BSO) \cong 
\zmod2$ and is generated by rank zero oriented vector bundles on 
$S^1 \times S^2$ and the non-bounding spin structure 
over the $S^1$ cartesian factor.  However, 
$\rttau (D_{\rho (\cdot )})$ is independent of the spin 
structure (as shown above) so that we can place the 
bounding spin structure on the $S^1$ cartesian factor to 
obtain the same element of $\T$.  But, that being done,
all of the structure bounds so that, in fact,
$\rttau_{S^1 \times S^2} (D_{\rho (\cdot )}) = 1$.
The cobordism invariance proves the proposition for a general 
compact 3-manifold.
\end{proof}

\subsection{spin-Chern-Simons and Chern-Simons}

We end our consideration of the classical Lagrangian field theory
by establishing its relation to Chern-Simons field theory; in particular,
as it was studied in \cite{F2}.   

\begin{proposition}\label{simply conn}
Let $\rho$ be a real rank zero virtual represenation of a connected, simply 
connected, compact Lie group $G$ and let $\pairing{}{}_{\rho}$ denote
the symmetric pairing defined by \eqref{p1 form}.  If $A$ is a 
$G$-connection over a closed, spin 3-manifold $X$, we let
$\exp 2 \pi i S_X (A)$ denote the $\T$-valued Chern-Simons 
invariant determined by the pairing $\frac{1}{2}\pairing{}{}_{\rho}$, 
as defined in \cite{F2}.  Then 
$$
\rttau_X (D_{\rho A}) = \exp 2 \pi i S_X (A).
$$ 
\end{proposition} 

\begin{proof}
The proof relies on the fact that the cobordism group $\Omega_3^{spin}(BG) = 0$
whenever $G$ is compact and simply connected.  In that case there exists a
spin 4-manifold $M$ such that $\d M = X$ as a spin manifold and there exists
an extension $A'$ of $A$ over $M$.  On the one hand it is well-know that
$$
\exp 2 \pi i S_X (\rho A) = 
\exp 2 \pi i \int_M \frac{1}{2}
\pairing{\Omega^{A'}}{\Omega^{A'}}_{\rho}.  
$$
On the other hand, since $\rho$ has rank zero the APS index theorem
implies
$$
\rttau_X (D_{\rho A}) =
\exp \pi i \int_M  
\pairing{\Omega^{A'}}{\Omega^{A'}}_{\rho},
$$
and this proves the proposition
\end{proof}

To make the equality of these theories even stronger we point out the
correspondence between their respective levels.
Recall that for Chern-Simons the levels are elements of $H^4(BG)$
while for spin-Chern-Simons the level are elements of $E^4(BG)$.
However, for $G$ simply connected there is a natural isomorphism
$i : H^4(BG) \rightarrow E^4(BG)$ so that Chern-Simons theory
at level $\alpha \in H^4(BG)$ is isomorphic to spin-Chern-Simons
at level $i (\alpha) \in E^4(BG)$.
              
\section{Spin-Chern-Simons on 3-Manifolds with Boundary}

In this chapter we extend our investigation of 
the Lagrangian spin-Chern-Simons field theory
to compact 3-manifolds with boundary.  When 
$X$ is without boundary our action $\rttau$
is a $\T$-valued function of the smooth parameters.
What we see in this chapter is that 
when $X$ has non-empty boundary $\d X$,
$\rttau$ is readily interpreted as a section
of the {\it Pfaffian line} associated to $\d X$. 
This interpretation and many of its consequences stem 
from previous work done for the $\tau$-invariant.
We review these less refined results before stating 
the corresponding results for $\rttau$.      

\subsection{A review of $\tau$ and DF-boundary conditions}

In this section we summarize some of the features of  
$\tau$-invariants for manifolds with boundary
as worked out in \cite{DF} and \cite{F3}.  
For the most part -- though much of this discussion
applies more generally in odd dimensions -- we focus on 
3-manifolds with boundary.

These results are formulated in terms of {\it graded} vector
spaces and especially graded lines.  For that reason we 
review a few salient points about graded lines 
from \cite{DF} and defer to the source for a more 
detailed account.   

The general situation we consider is a line of the form
$\Det V = \bigwedge^{\dim V} V$ for some 
finite dimensional complex vector space $V$ and we assign 
to it the grading $|\Det V| = \dim V$ (mod 2).  The 
gradings add or substract when we tensor two lines
or tensor one by the inverse of the other, respectively.  
One of the ways in which gradings affect our computations 
is when we consider the natural isomorphisms :
$$
L^* \otimes L \rightarrow \C
\quad \text{and} \quad
L \otimes L^* \rightarrow \C 
$$
where $L$ is a graded line with grading $|L| \in \zmod2$. 
In the graded category we take the convention that the
first isomorphism $a^{-1} \otimes b \mapsto a^{-1}(b)$
does not involve any signs.  However, this implies that
the second isomorphism must take the form $b \otimes 
a^{-1} \mapsto (-1)^{|L|} a^{-1}(b)$.  The second
isomorphism is called the {\it supertrace} and is denoted
$\sTr$.  Notice that if $L$ has an even grading it is 
equivalent to the normal trace and if $L$ has an odd 
grading then it is equivalent to minus the normal trace.

Let $X$ be a spin Riemannian 3-manifold with boundary
$\d X$.  Recall that $\tau$ is a spectral invariant of 
the (possibly twisted) Dirac operator.  When $\d X = \emptyset$
the Dirac operator $D_X$ is self-adjoint and elliptic so that
$\tau$ is well-defined.  To maintain these analytic properties
when $\d X \neq \emptyset$ Dai and Freed impose 
elliptic boundary conditions similar to those introduced by 
Atiyah-Patodi-Singer but require an additional piece of data to
adapt the APS boundary conditions to odd-dimensional manifolds. 
We explain this next.
 
Since $\d X$ is even-dimensional its spinor bundle 
$S_{\d X}$ split as $S_{\d X} = S^+_{\d X}
\oplus S^-_{\d X}$ and the Dirac operator $D_{\d X} :
C^{\infty}(S^{\pm}_{\d X}) \rightarrow C^{\infty}
(S^{\mp}_{\d X})$ interchanges the components.  The cobordism 
invariance of the index implies that $\dim \text{Ker}^+ 
D_{\d X} = \dim \text{Ker}^- D_{\d X}$.  The additional piece
of data is an isometry
$$
T : \text{Ker}^+ D_{\d X} \longrightarrow
\text{Ker}^- D_{\d X}.
$$ 
These boundary conditions then ensure that $\tau (D_X)$ is
well-defined.  Though it does depend on the isometry $T$ its 
dependence can be factored out in such a way that
$$
\tau_X \in \Det_{\d X}^{-1}
$$
where $\Det_{\d X}$ is the {\it determinant line}
of the Dirac operator $D_{\d X}$:
\begin{equation}\label{det line}
\Det_{\d X} = (\Det \text{Ker}^- D_{\d X} ) \otimes
(\Det \text{Ker}^+ D_{\d X} )^{-1} .
\end{equation}
Here and in what follows  $L^{-1} = L^*$ 
for any (abstract) line $L$.  Furthermore, $| \tau_X |^2 = 1$ in
the {\it Quillen metric} on $\Det_{\d X}^{-1}$.  
From the definition \eqref{det line} we see that the determinant 
line $\Det_{Y}$ associated to a closed, even dimensional Riemannian 
spin manifold $Y$ has the grading index$(D_Y)$ (mod 2).  

\begin{remark}\label{virt detline}
Above we (parenthetically) stated that all of this was 
possible in the context of Dirac operators twisted 
by virtual vector bundles.  For clarity, and becuase
it is essential to our considerations below, we make
explicit what we mean.  First we point out that if
$(E_j, \nabla^{E_j})$ $j = 1,2$ are two vector bundles 
with unitary connections over a spin Riemannian 
3-manifold $X$ then \eqref{det line} implies 
$$
\Det_{\d X, E_1 \oplus E_2} = \Det_{\d X, E_1} 
\otimes \Det_{\d X, E_2}.
$$
The definition that is compatible with addition of virtual
vector bundles is therefore
$$
\Det_{\d X, E_1 - E_2} = \Det_{\d X, E_1} 
\otimes \Det^{-1}_{\d X, E_2}
$$
and so $\tau_X(D_{E_1 - E_2})$ is an element of 
of $\Det^{-1}_{\d X, E_1 - E_2}$. 
\end{remark}

If $X \rightarrow Z$ is a family of Riemannian 3-manifolds
with boundary then $\d X \rightarrow Z$ is a family of 
2-manifolds.  (Also allowed is a twisting virtual vector bundle 
$E \rightarrow X$with unitary connection $\nabla^E$).  
The lines \eqref{det line} patch together to form a smooth 
line bundle $\Det_{\d X/Z} \rightarrow Z$ with the Quillen 
metric and a natural unitary connection $\nabla'$, defined in 
\cite{BF1}.  In this context we have that
$$
\tau_{X/Z} : Z \longrightarrow \Det^{-1}_{\d X/Z} 
$$
is a smooth unitary section.  One of the basic results 
regarding this section is a variation formula which computes
its covariant derivative over the family.
\begin{theorem}[\cite{DF}, Theorem 1.9]
With respect to the natural connection $\nabla'$ on 
$\Det^{-1}_{\d X/Z}$,
\begin{equation}\label{var tau}
\nabla' \tau_{X/Z} = 2 \pi i 
\left[ \int_{X/Z} \widehat{A}(\Omega^{X/Z}) ch (\Omega^E)  
\right]_{(1)} \tau_{X/Z}.
\end{equation}
\end{theorem} 
Notice that this is just the (exponentiated) generalization of 
\eqref{difofxi}.  

For now, this completes our review of the results
in \cite{DF}.  We will continue this discussion
later in this chapter when we review a gluing law
for $\tau_X$ and then again in the next chapter when
we consider Hamiltonian spin-Chern-Simons theory over 
closed 2-manifolds.

\subsection{$\rttau$ and DF-boundary conditions}

We come back to $\rttau$ which is the spectral invariant 
relevant to our spin-Chern-Simons field theory.  With some care but
little difficulty the results summarized above for $\tau$
can be refined for $\rttau$.  We explain this next.  Though
much of what we say generalizes to 3 (mod 8) dimensional 
manifolds (thanks to Bott periodicity), we maintain our focus 
on 3-manifolds with boundary.

Let $X$, once again, be a spin Riemannian 3-manifold with 
boundary $\d X$.  Recall that when $\d X = \emptyset$
the spectral invariant $\rttau$ is well defined because 
the (possibly twisted) Dirac operator is self-adjoint, elliptic,
and compatible with a quaternionic structure $J$ on the 
the space $C^{\infty}(S_X)$ of spinor fields on $X$.  
To maintain these analytic properties when $\d X \neq 
\emptyset$, including quaternionic compatibility, we 
require a refinement of the DF boundary conditions.  
We explain this next.

Since $\d X$ is a closed 2-manifold the spinor bundles
$S^+_{\d X}$ and $S^-_{\d X}$ are naturally dual to 
each other and the Dirac operator $D_{\d X} : C^{\infty}
(S^+_{\d X}) \rightarrow C^{\infty}(S^-_{\d X})$ is 
complex skew-symmetric with respect to that duality.  This 
informs our requirements on the isometry
$$
T : \text{Ker}^+D_{\d X} \longrightarrow 
\text{Ker}^-D_{\d X}.
$$
The cobordism invariance of the mod-2 index implies that 
$\dim \text{Ker}^+D_{\d X}$ is even dimensional and so
supports a skew-symmetric isometry from itself to its dual.
We require that the isometry $T$ be skew-symmetric.  

These boundary conditions ensure that $\rttau (D_X)$ is
well-defined.  Indeed, restricted the subspace of spinors
that satisfy the boundary conditions, $D_X$ is elliptic,
self-adjoint, and quaternionic.  Though it depends on $T$ that 
dependence can be factored out so that one observes
$$
\rttau_X \in \Pfaff^{-1}_{\d X}
$$
where $\Pfaff_{\d X}$ is the {\it Pfaffian line} of the 
Dirac operator $D_{\d X}$:
\begin{equation}\label{pf line}
\Pfaff_{\d X} = \Det \text{Ker}^- D_{\d X}
\end{equation} 
Notice that, because of the duality mentioned above,
there is a natural identification 
\begin{equation}\label{sqpf=det}
\Pfaff^{\otimes 2}_{\d X} = \Det_{\d X}.
\end{equation}  
Under this identification we have (not surprisingly) 
$\rttau_X \otimes \rttau_X = \tau_X$
This implies $| \rttau_X |^2 = 1$ in the
square root of the Quillen metric on $\Pfaff^{-1}_{\d X}$.
We also point out that the Pfaffian line 
$\Pfaff_{Y}$ associated to a closed, 2 (mod 8)-dimensional 
Riemannian spin manifold $Y$ has the grading $\text{ind}_2 
(D_Y) = \dim \text{Ker}^+ D_Y$ (mod 2).

\begin{remark}\label{virt pfline}
Remarks analogous to those made in Remark \ref{virt detline} 
apply to virtual vector bundles and Pfaffian lines.
Indeed, if $(E_j, \nabla^{E_j})$ $j = 1,2$ are two real vector 
bundles with orthogonal connections over a spin Riemannian 
3-manifold $X$ then \eqref{pf line} implies 
$$
\Pfaff_{\d X, E_1 \oplus E_2} = \Pfaff_{\d X, E_1} 
\otimes \Pfaff_{\d X, E_2}.
$$
The definition that is compatible with addition of virtual
vector bundles is therefore
$$
\Pfaff_{\d X, E_1 - E_2} = \Pfaff_{\d X, E_1} 
\otimes \Pfaff^{-1}_{\d X, E_2}.
$$
and so $\rttau_X(D_{E_1 - E_2})$ is an element of 
of $\Pfaff^{-1}_{\d X, E_1 - E_2}$.  
\end{remark}

Let $X \rightarrow Z$ be a family of Riemannian 3-manifolds 
with boundary and $\d X \rightarrow Z$ the corresponding
family of closed 2-manifolds.  (Also allowed is a twisting real
virtual vector bundle $E \rightarrow X$ with orthogonal 
connection $\nabla^E$).  The lines \eqref{pf line} patch
together to form a smooth line bundle $\Pfaff_{\d X/Z} 
\rightarrow Z$ with the square root of the Quillen metric 
and a natural unitary connection $\nabla$ \cite{F1}.  In 
this context we have that 
$$
\rttau_{X/Z} : Z \longrightarrow \Pfaff^{-1}_{\d X/Z}
$$ 
is a smooth unitary section.  The identification \eqref{sqpf=det}
also patches together so that $\Pfaff^{\otimes 2}_{\d X/Z} =
\Det_{\d X/Z}$ and $\rttau_{X/Z} \otimes \rttau_{X/Z}
= \tau_{X/Z}$.  Furthermore, we can identify the connections
$\nabla^{\otimes 2} = \nabla'$.  This implies a variation formula
for the section $\rttau_{X/Z}$ which follows trivially from 
\eqref{var tau}.
\begin{proposition}
With respect to the natural connection $\nabla$ on 
$\Pfaff^{-1}_{\d X/Z}$,
\begin{equation}\label{var rttau}
\nabla \rttau_{X/Z} =  \pi i 
\left[ \int_{X/Z} \widehat{A}(\Omega^{X/Z}) ch (\Omega^E)  
\right]_{(1)} \rttau_{X/Z}.
\end{equation}
\end{proposition}

Actually, many of these considerations are too general for our
needs (they are true in all dimensions 3 (mod 8) for instance)
and so we specialize to spin-Chern-Simons Lagrangian
field theory next.

\subsection{Functoriality of $\rttau_X \in \Pfaff^{-1}_{\d X}$}

Let $G$ be a compact Lie group and $\rho $ a real oriented virtual
representation of $G$.
Let $Y$ be a closed spin Riemannian 2-manifold with Dirac operator
$D_Y : \C^{\infty}(S^+) \rightarrow \C^{\infty}(S^-)$.  We define 
$\Cat{G}{Y}$ to be the category of $G$-connections over $Y$.  We 
further define a functor from $\Cat{G}{Y}$ to the category of complex 
lines $\mathscr{L}$
\begin{align}
\Line{\rho}{} : \Cat{G}{Y} & 
\longrightarrow \mathscr{L} \nonumber \\
(Q,B) & \longmapsto \Pfaff^{-1}_Y(D_{\rho B})
\label{line functor}
\end{align} 
where $\Pfaff_Y^{-1}(D_{\rho B})$ is the inverse Pfaffian line
of the twisted Dirac operator $D_{\rho B}$ on $Y$, as defined in 
\eqref{pf line}.  From the discussion above, if $B_u$ is a smooth
family of $G$-connections on $Y$ varying over a smooth manifold
$U$ then the inverse Pfaffian lines $\Line{\rho}{(B_u)}$ form 
a smooth hermitian line bundle with unitary connection over $U$.

We come back to a spin Riemannian 3-manifold $X$ with boundary
$\d X$ and Dirac operator $D_X : \C^{\infty}(S) \rightarrow 
\C^{\infty}(S)$.  From the discussion above, if $\d X \neq \emptyset$,
and $\Cat{G}{X,B}$ is the category of $G$-connections $A$ on $X$ 
such that $A|{\d X} = B$, then we have the assignment
\begin{align*}
\Cat{G}{X, B} & \longrightarrow \Line{\rho}{(B)} \\
A & \longmapsto \rttau_X (D_{\rho A})
\end{align*}
From the discussion above if $A_u$ is a smooth family of 
$G$-connections on $X$ varying over a smooth manifold $U$ then
the elements $\rttau_X (D_{\rho A_u}) \in \Line{\rho}
{(A_u|_{\d X})}$ form a smooth section of line bundle formed
by the lines $\Line{\rho}{(A_u|_{\d X})}$.  To allow for 
the case when $X$ is closed we set $\Line{\rho}{} = \C$ when 
$Y = \emptyset$.  If $X$ itself is empty then we set 
$\rttau_X = 1$.  In light of this new interpretation for
the action of our field theory, we present a properly 
adjusted version of Proposition \ref{functorial}.

\begin{proposition}\label{functorial2}
Let $G$ be a compact Lie group and $\rho$ a virtual 
real representation of $G$.  Then the assignments
\begin{align*}
B \longmapsto \Line{\rho}(B), & 
\quad B \in \Cat{G}{Y} \\
A \longmapsto \rttau_X (D_{\rho A}), &
\quad A \in \Cat{G}{X},
\end{align*}
defined above for closed spin Riemannian 2-manifolds $Y$
and compact spin Riemannian 3-manifolds $X$ are smooth
and satisfy:
\begin{enumerate}

\item (Functoriality)If $\psi: Q' \rightarrow Q$ and 
$H: S'_{g', \sigma'} \rightarrow S_{g, \sigma }$ are any 
G-bundle and spinor bundle morphisms (respectively) covering an 
orientation and spin structure preserving isometry $h: (Y',g') 
\rightarrow (Y,g)$, and $B$ is a connection on $Q$, then there
is an induced isometry 
\begin{equation}\label{isom on L}
(H \otimes \rho \psi)^* : \Line{\rho}(B)
\longrightarrow \Line{\rho}{(\psi^*B)}
\end{equation} 
and these compose properly.  If $\varphi: P \rightarrow P$ and 
$F: S'_{g', \sigma'} \rightarrow S_{g, \sigma }$ are any 
G-bundle and spinor bundle morphisms (respectively) covering an 
orientation and spin structure preserving isometry $f: (X',g') 
\rightarrow (X,g)$, and $A$ is a connection on P, then
\begin{equation}\label{isom on rttau}
\d (F \otimes \rho \varphi )^* \rttau_{X} (D_{\rho A}) =
 \rttau_{X'} (D'_{\rho (\phi^*A)})
\end{equation}
where $D'$ is the Dirac operator on $X'$ and 
$\d (F \otimes \rho \varphi ) : \d (S_{X'} \otimes \rho P) 
\rightarrow  \d (S_{X} \otimes \rho P)$ is the induced map over
the boundary. 
\item (Orientation) There is a natural isometry 
\begin{equation}\label{line orient}
\Line{\rho}{(-Y,B)} \cong {\Line{-\rho}{(Y,B)}}
\end{equation} 
such that
\begin{equation}\label{rttau orient}
\rttau_{-X}(D_{\rho A}) = (-1)^{\binom{k}{2}}
\rttau_X(D_{\rho A})^{-1}.
\end{equation}
where $k$ is the number of components $Y \subset \d X$ such that
the twisted (chiral) Dirac operator on $Y$ has non-trivial mod 2 index. 
\item (``Additivity'')If $Y = Y_1 \sqcup Y_2$ is a disjoint 
union, and $B_j$ are connections over $Y_j$, then there is 
a natural isometry
\begin{equation}\label{line disj}
\Line{\rho}{(B_1 \sqcup B_2)} \cong 
\Line{\rho}{(B_1)} \otimes \Line{\rho}{(B_2)}. 
\end{equation} 
If $X = X_1 \sqcup X_2$ is a disjoint union,
and $A_j$ are connections over $X_j$, then
$$ 
\rttau_{X_1 \sqcup X_2}(A_1 \sqcup A_2) =
\rttau_{X_1}(A_1) \otimes \rttau_{X_2}(A_2).  
$$
\end{enumerate}  
\end{proposition}
Here and in what follows if $a \in L$ is an 
element of a line then we denote by $a^{-1}
\in L^{-1}$ the unique dual element such that
$a^{-1}(a) = 1$.  In particular, this applies
to part (2) of the proposition.

\begin{remark}\label{func impl}
Let us pause to take stock of what is implied 
by the functorial statements in this proposition.
When we defined the assignment 
$$
\Line{\rho}{}: \Cat{G}{Y} \longrightarrow \Line{}{}
$$
we stated, without proof, that it was a functor.
Part (1) of the proposition proves this claim.
From a less lofty point of view, if $Q \rightarrow Y$
is a $G$-bundle over a closed spin Riemannian 2-manifold
there is an associated line bundle $\Line{\rho}{(Q)}
\rightarrow \Cat{}{Q}$ over the space of connections,
and the action of the gauge transformations $\Mor{}{Q}$
lifts to $\Line{\rho}{(Q)}$.  Also, a $G$-bundle
$P \rightarrow X$ over a compact spin Riemannian 
3-manifold determines a restriction map $\Cat{}{P} 
\rightarrow \Cat{}{\d P}$ and so a pulled back line
bundle $\Line{\rho}{(P)} \rightarrow \Cat{}{P}$.  The
action of the gauge transformations $\Mor{}{P}$ lift
to $\Line{\rho}{(P)}$ and $\rttau_X$ is an invariant
section of $\Line{\rho}{(P)}$. We will have more to
say about the implications of this proposition
once we have stated the gluing law \eqref{glue rttau}
for $\rttau_X$.
\end{remark}

\begin{proof}
To prove (1) let $\Phi_{H,\psi}$ denote
the isometry on sections induced by $H \otimes \rho 
\psi$.  Then $\Phi_{H,\psi}$ restricts to an isometry
$\text{Ker}^- D_{\rho B} \rightarrow 
\text{Ker}^- D'_{\rho (\psi^*B)}$, where $D'$ is the 
Dirac operator on $Y'$.  We denote the restriction
by  $\Phi_{H,\psi}$ as well so that $\Det \Phi_{H,\psi}$
is the isometry in \eqref{isom on rttau}.  That $\rttau$ 
obeys \eqref{isom on rttau} follows from the proof for
part (1) of Proposition \ref{functorial} and the fact
that the induced isometry on sections $\Phi_{F,\varphi}$
clearly preserves the (refined) DF-boundary conditions.

To prove (2) recall that under change in orientation
there is a switch in chirality, $S_{-Y}^{\pm} = S_Y^{\mp}$.
The duality of $S_Y^{+}$ and $S_Y^{-}$, together 
with the comments in Remark \ref{virt pfline}, imply 
\eqref{line orient}.  That $\rttau$ obeys 
\eqref{rttau orient} follows partly from the proof for part
(2) of Proposition \ref{functorial} and the fact that
the DF boundary conditions behave appropriately upon change
in orientation.  For the sign factor, recall that the isomorphism is 
one between graded lines and so will involve signs.  In particular,
identifying $\Line{-1}{}_1 \gotimes \dots \gotimes 
\Line{-1}{}_n$ with $(\Line{-1}{}_1 \gotimes \dots 
\Line{-1}{}_n)$ involves permuting $\binom{k}{2}$
odd elements past each other whereby we obtain the sign
in \eqref{rttau orient}.   

The proof for (3) is simple and left to the reader.    
\end{proof}
   
We offer a properly adjusted version of Proposition 
\ref{ext of grp} for 3-manifolds with boundary.

\begin{proposition}
Let $i : G \hookrightarrow G'$ be an inclusion of compact 
Lie groups.  Suppose $\rho$ is real virtual representation of 
$G'$.  Then if $(Q,B)$ is a $G$-connection over a closed, spin 
Riemannian 2-manifold $Y$, and $(Q',B')$ its $G'$ extension, 
there is a natural isometry
\begin{equation}\label{incl isom}
i_B : \Line{\rho \circ i}{(B)} \longrightarrow
\Line{\rho}{(B')}.
\end{equation}
If $A$ is a $G$-connection over a compact, spin Riemannian 
3-manifold $X$, and $A'$ its $G'$ extension, then
$$
i_B \left( \rttau_X (D_{(\rho \circ i ) B}) \right) 
= \rttau_X (D_{\rho B'})
$$
\end{proposition} 

In categorical language, $i$ induces a tranformation from the
functor $\Cat{G}{Y} \rightarrow \Line{}{}$ to the functor
$\Cat{G}{Y} \rightarrow \Cat{G'}{Y} \rightarrow 
\Line{}{}$.  For each $X$ this induces a transformation from 
$\Cat{G}{X} \rightarrow \Cat{G}{\d X} \rightarrow 
\Line{}{}$ to $\Cat{G}{X} \rightarrow \Cat{G'}{X} 
\rightarrow \Cat{G'}{\d X} \rightarrow \Line{}{}$,
that preserves the elements $\rttau_X$.  

\begin{proof}
There is a natural isomorphism between the associated virtual 
vector bundles $\rho Q' \rightarrow (\rho \circ i)Q$ which
sends the associated virtual connection $\rho B'$ to $(\rho 
\circ i) B$.  The induced isometry between the sections 
restricts from the kernel of $D_{\rho B'}$ to the kernel 
$D_{(\rho \circ i)B}$.  If we let $\Phi$ denotes the restriction
then $\Det \Phi$ provides the isometry \eqref{incl isom}.
The rest of the proposition follows from the proof Proposition
\ref{ext of grp}.   
\end{proof}

\subsection{Gluing Formulae}

We return to the results of \cite{DF}.  In particular,
we review the gluing formula for the $\tau$-invariant.
The context of the gluing formula is as follows.  Let $X$
be a compact, spin Riemannian 3-manifold and $Y 
\hookrightarrow X$ a closed, spin 2-dimensional 
submanifold.  We cut along $Y$ to obtain a new manifold
$X^{\text{cut}}$ with $\d X^{\text{cut}} = \d X 
\sqcup Y \sqcup -Y$.  Recall that 
$$
\tau_X \in \Det^{-1}_{\d X}
\quad \text{and} \quad
\tau_{X^{\text{cut}}} \in \Det^{-1}_{\d X} 
\otimes \Det^{-1}_{Y} \otimes \Det^{-1}_{-Y}
$$ 
and notice that $\Det^{-1}_{-Y} \cong (\Det^{-1}_{Y})^{-1}$.
Thus we can apply the supertrace to the last two factors 
of the tensor product to obtain an element of 
$\Det^{-1}_{\d X}$.  That being said we state the gluing
formula for $\tau$.

\begin{theorem}[\cite{DF}, Theorem 2.20]
In the context described above
\begin{equation}\label{glue tau}
\sTr (\tau_{X^{\text{cut}}}) = \tau_X .
\end{equation}
\end{theorem}

With some care but very little difficulty we can 
refine the discussion and results just given to Pfaffian lines 
and the $\rttau$-invariant.  Recall that the Pfaffian line 
$\Pfaff_{Y}$ associated to a closed, 2 (mod 8)-dimensional 
Riemannian spin manifold $Y$ has the grading $\text{ind}_2 
(D_Y)$.  With that in mind, in the context of the gluing formula 
we have that
$$
\rttau_X \in \Pfaff^{-1}_{\d X}
\quad \text{and} \quad
\rttau_{X^{\text{cut}}} \in \Pfaff^{-1}_{\d X} 
\otimes \Pfaff^{-1}_{Y} \otimes \Pfaff^{-1}_{-Y}.
$$ 
We apply the supertrace to the last two factors of the tensor
prodect to obtain the gluing formula for $\rttau$.
\begin{proposition}
In this case
\begin{equation}\label{glue rttau}
\sTr (\rttau_{X^{\text{cut}}}) = \rttau_X .
\end{equation}
\end{proposition}

This very general gluing formula, of course, applies 
to the action of our spin-Chern-Simons theory.  Due
to its importance in the upcoming chapter on the Hamiltonian
field theory we restate the gluing formula for $\rttau$ in the 
notation of the previous section.  The reader should consider
it as an addendum to Proposition \ref{functorial2}.

\begin{corollary}[Gluing]
Let $X$, $Y$, $G$, and $\rho$ be as in Proposition 
\ref{functorial2}.  Now suppose $Y \hookrightarrow X$ 
is a closed spin submanifold and $X^{\text{cut}}$ is the 
manifold obtained by cutting $X$ along $Y$.  Then 
$\d X^{\text{cut}} = \d X \sqcup Y \sqcup -Y$.  Suppose
$A$ is a $G$-connection over $X$, with $A^{\text{cut}}$ the
induced connection over $X^{\text{cut}}$, and $B = A|_Y$.
Then 
\begin{equation}\label{glue}
\rttau_X (D_{\rho A} ) =
\sTr  \left( \rttau_{X^{\text{cut}}}
( D_{\rho A^{\text{cut}}} ) \right)
\end{equation}
where $\sTr $ is the contraction
$$
\sTr : \Line{\rho}{(A^{\text{cut}})}
\cong \Line{\rho}{(\d A)} \otimes \Line{\rho}{(B)}
\otimes \Line{-\rho}{(B)} \longrightarrow
\Line{\rho}{(\d A)}
$$
taking the supertrace of the last two factors of the
tensor product.
\end{corollary}

\begin{remark}
We return to the comments in Remark \ref{func impl}.  
Part (3) of Proposition \ref{functorial2} and the gluing 
formula expresses the fact that the action is a local functional of local
fields.  Part (1) implies that the action is invariant under the 
symmetries of the fields, and part (2) expresses the fact that the
action is unitary.    
\end{remark}

\subsection{Dependence of $\rttau$ on smooth parameters: redux}

Let $X$ be a compact, spin 3-manifold with a possibly non-empty 
boundary $\d X$.  In Propositions \ref{d wrt metric} and 
\ref{d wrt conn} we compute how the {\it function} $\rttau$ 
varied with respect to infinitesimal changes in the metric and 
connection.  This was done assuming $\d X = \emptyset$.  
We now offer analogous propositions, adjusted to the case $\d X \neq 
\emptyset$, so that we compute how the {\it section} $\rttau$     
varies with respect to infinitesimal changes.  Of course, we make this
computation using the natural connection on the Pfaffian line bundle, 
whose properties we reviewed above.  In particular, we apply the 
variation formula \eqref{var rttau} to the section
determined by the action of our spin-Chern-Simons theory.

Let $P \rightarrow X$ be any principal $G$-bundle.  
Then we take the obvious family of Riemannian 3-manifolds with 
$G$-connections
$$
\Met (X) \times \Cat{}{P} \times X 
\longrightarrow \Met (X) \times \Cat{}{P} = Z.
$$ 
Given a real virtual representation $\rho$ 
we obtain a family of twisted Dirac operators parametrized
by $Z$.  From Remark \ref{func impl}  
we know that there is a line bundle $\Line{\rho}{(X, P)} 
\rightarrow Z$ such that the action is a unitary section
$\rttau : Z \rightarrow \Line{\rho}{(X, P)}$.
Recall that $\Line{\rho}{(X,P)}$ is a bundle of 
Pfaffian lines and therefore has a natural connection $\nabla$. 
Having said that, we can state how the action behaves 
under infinitesimal variations.  We begin with a generalization 
of Proposition \ref{d wrt metric}.
  
\begin{proposition}\label{d wrt metric2}
If $\rho $ is a real rank zero virtual representation  
then, with respect to $\nabla$,  $\rttau$ is covariantly 
constant along $\Met (X)$. 
\end{proposition}

\begin{proof}
Using the variation formula \eqref{var rttau} the  
proof of Proposition \ref{d wrt metric} carries over almost word
for word to this situation.
\end{proof}

Now we generalize Proposition \ref{d wrt conn}.

\begin{proposition}\label{d wrt conn2}
Let $A_t$ be a path in $\Cat{}{P}$ and let $\Omega_t$ denote
the curvature of the connection $A_t$, $\dot{A_t}$ the tangent
vector along the path.  Then
$$
\nabla_{\dot{A_t}} \rttau (A_t) =
2 \pi i \cdot dt \cdot \rttau (A_t)
\int_{X} \fpairing{\Omega_t}{\dot{A_t}}_{\rho}.
$$
\end{proposition}

\begin{proof}
Once again, using the variation formula \eqref{var rttau}
the proof of Proposition \ref{d wrt conn} carries carries over
almost word for word to this situation.
\end{proof}

\begin{remark}\label{mod space1}
Assuming $\pairing{}{}_{\rho}$ is nondegenerate, this 
proposition implies that $\nabla \rttau |_A = 0$ if and only if
$\Omega^A = 0$; that is, if and only if $A$ is flat.  
Since the section $\rttau$ is invariant with
respect to gauge transformations $\Mor{G}{X}$, so too is
the space of solutions to the Euler-Lagrange equation.
This is also obvious from \eqref{eul-lag}.  We let 
$$
\Mod{G}{X} \subset \overline{\Cat{G}{X}}
$$
denote the space of equivalence classes of solutions
to \eqref{eul-lag}.  We will have more to say about
$\Mod{G}{X}$ when we discuss the Hamiltonian
field theory in Section 1.3. 
\end{remark}

\section{The Classical Hamiltonian Theory}

In the Lagrangian spin-Chern-Simons theory the ``spacetimes''
are compact, spin, Riemannian 3-manifolds, possibly with
non-empty boundary.  In the Hamiltonian field theory the spacetimes 
are spin, Riemannian 3-manifolds that are globally products of
a closed, spin 2-manifold (``space'') and an infinite interval 
(``time'').  The classical solutions to the Euler-Lagrange equations 
are flat connections and we are only interested in equivalence
classes of these with respect to the gauge group.  This is the
classical phase space.  In the Hamiltonian theory the solutions 
are constant with respect to time so that the classical phase space
is the moduli space of flat connections on the 2-manifold.  
These spaces have been the topic of much study in the past 
twenty years.  Since our action is defined using certain elements 
of geometric index theory, it leads quite naturally to 
appearance of Pfaffian line bundles with their natural geometry:
a metric and a unitary connection.  We have already
seen foreshadowings of this in the previous section.  Here the
line bundles themselves will take center stage, and their 
geometry will inform our formulation of the classical
Hamiltonian field theory.

\subsection{Moduli spaces of flat connections}

Recall that in the Lagrangian field theory over a compact
spin 3-manifold $X$, the classical space of solutions to the 
Euler-Lagrange equations is the moduli space of flat connections
$\Mod{G}{X}$.  This space appears again in the Hamiltonian
field theory.  For that reason we review some standard facts about 
moduli spaces of flat connections.  They reflect on how
the classical theory (both Lagrangian and Hamiltonian)
probes the topology of the spacetimes.  

\begin{proposition}
Let $X$ be any smooth manifold and let $\{ x_i \}_{i \in \pi_0(X)}$
be a set of basepoints for each component of $X$.  Then the holonomy
provides a natural identification
\begin{equation}\label{mod space}
\Mod{G}{X} = \prod_{i} \Hom (\pi_1(X,x_i),G)/G,
\end{equation}
where $G$ acts on $\Hom (\pi_1(X,x_i),G)$ by conjugation.  
Furthermore, this identification is independent of the basepoints.
\end{proposition}
The proof is standard.

Typically the moduli space $\Mod{G}{X}$ is not a manifold.
However, if $Y$ is a compact, oriented  2-manifold then it is
well known that $\Mod{G}{Y}$ is a stratified space and that 
the stratum of top dimension is a smooth manifold \cite{G},
\cite{AB}.  To investigate this manifold structure near the equivalence 
class of a flat connection $A$, we consider the twisted complex
\begin{equation}\label{A complex}
0 \rightarrow \Omega_X^0(\text{ad}P) 
\overset{d_A}{\longrightarrow}
\Omega_X^1(\text{ad}P) \rightarrow \dots .
\end{equation}
Here $d_A$ is the usual extension of $A$ to 
act on differential forms with value in $\text{ad}P$; and that this 
is a complex follows from $d_A^2 = \Omega^A = 0$.
If we denote the cohomology groups of this complex by 
$H^{\bullet}(X; d_A)$, then one might guess that the tangent
space at $A$ is 
\begin{equation}
T_A \Mod{G}{X} \cong H^1(X, d_A).
\end{equation}
Indeed, this is certainly the case whenever $A$ represents a smooth
point of $\Mod{G}{X}$.

Many of the properties of the usual deRham cohomology carry over
for twisted cohomology.  For example, if $Y$ is a compact, oriented
2-manifold then for a flat connection $B$ there is a nondegenerate
pairing 
$$
H^0(Y,d_B) \otimes H^2(Y,d_B) \longrightarrow \R ,
$$
which mimics Poincare duality in the usual deRham cohomology.
The zeroth cohomology $H^0(Y,d_B)$ is the Lie algebra of the 
stabilizer $\Mor{}{B}$.  Of course, it vanishes if $B$ is irreducible
since, in that case, the stabilizer of $B$ is the center of $G$ which 
is finite.  At smooth points the index theorem for the twisted
complex gives
\begin{equation}\label{dim formula}
\dim \Mod{G}{Y} = \dim H^1(Y,d_B) = 
- \dim G \cdot \chi (Y) + 2 \dim \Mor{}{B}. 
\end{equation}
When $Y$ has a complex structure the flat $G$-connections can identified
with holomorphic structures on $G_{\C}$ principal bundles and this imbues 
$\Mod{G}{Y}$ with its own complex structure \cite{NS}, \cite{AB}.  
This is manifest in the dimension \eqref{dim formula} which is always 
even.

\begin{remark}\label{stacks1}
More properly, what we actually work with is the moduli {\it stack} 
of flat connections.  We briefly explain what we mean.  Recall that the 
space of fields is the category of $G$-connections $\Cat{G}{X}$.  
The objects are pairs $(P,A)$ -- where $P \rightarrow X$ is a $G$-bundle
and $A$ is a connection on $P$ -- and the morphisms $\Mor{G}{X}$
are $G$-bundle isomorphisms that cover the identity of $X$.
As all of the elements of $\Mor{G}{X}$ are invertible, the
category $\Cat{G}{X}$ is a {\it groupoid}.  The subcategory of flat
$G$-connections is a {\it subgroupoid}.
In working with the moduli stack, instead of the moduli space,
we keep track of, not just the equivalence class of $A$, but the 
automorphisms of $A$ as well.  That is, if $A$ has a non-trivial 
stabilizing subgroup $\Mor{}{A} \subset \Mor{}{P}$, then the 
moduli stack keeps track of that information; whereas the moduli 
space only sees the equivalence class of $A$.  Notice that if we 
only work with those connections with trivial stabilizing subgroups
then working with the moduli stack provides the same information
as working with the moduli space.  
\end{remark}  

\subsection{The space of fields for the Hamiltonian theory}

We consider a 3-manifold with boundary of the form
$X = [0,\infty) \times Y$ where $Y$ is a closed spin 
Riemannian 2-manifold.  $X$ has a natural product spin
structure determined entirely by the spin structure on $Y$,
and we impose the cylindrical metric (flat in the $[0,\infty)$
direction) on $X$.  The Dirac operator on $X$, even with
the DF-boundary conditions, is not an elliptic operator.  
In particular, the spectrum is not entirely discrete so that 
the $\rttau$-invariant is not defined.  Nonetheless, we 
can consider the critical points by taking compactly supported
variations.  Indeed, the log derivative $d \rttau_X / \rttau_X$ 
has an expression \eqref{d wrt conn} in terms of local fields that 
is well defined when the variation is taken along compactly supported 
fields.  Thus, we define a critical $G$-connection $A$ to be one for which
$$
\Omega^A = 0.
$$  
In the Hamiltonian formulation we may interpret the space 
of fields as a space of paths.  To do so we require the following proposition.

\begin{proposition}[\cite{F2} Proposition 3.14] 
Let $\{ Q \}$ be a set of representations for the
equivalence classes of principal $G$ bundles over $Y$.
Then there is an identification
\begin{equation}\label{path}
\overline{\Cat{G}{[0,\infty) \times Y}} =
\bigsqcup_{\{ Q \} } 
\Map ([0,\infty), \Cat{}{Q}) / \Mor{}{Q}.
\end{equation}
\end{proposition}

The proof makes use of the fact that $[0,\infty)$ is contractible
so that the topological type of a G-bundle $P \rightarrow  [0,\infty) 
\times Y$ is determined by $\d P \rightarrow Y$.  It also uses the 
fact that every connection on $ [0,\infty) \times Y$
is isomorphic to a connection that is trivial in the $[0,\infty)$ 
direction, which is exactly a path of connections on $Y$.  
In the physics literature, this is the statement that we can always work
in a ``temporal gauge'' in which the connection is trivial along time direction.  

If $A$ is a flat connection then $\d A = A|_Y$ is also flat.
The next proposition implies that the equivalence class of the 
classical solutions on $ [0,\infty ) \times Y$ are completely 
determined by the equivalence class of their initial value.

\begin{proposition}[\cite{F2} Proposition 3.16]
The restriction to the boundary 
$$
\Mod{G}{[0,\infty) \times Y}
\subset \overline{\Cat{G}{[0,\infty) \times Y}}
\longrightarrow \overline{\Cat{G}{ \{ 0 \} \times Y}}
$$
is an isomorphism of $\Mod{G}{[0,\infty) \times Y}$ onto
the moduli space $\Mod{G}{Y}$ of flat connections over $Y$.
\end{proposition} 

In Hamiltonian classical mechanics the fields over a cylinder
are paths in a {\it symplectic} manifold.  We address the 
issue of symplectic structure next by bringing the Pfaffian lines 
$\Line{\rho}{}$ to the forefront. 

\subsection{Geometry of the Pfaffian line bundle}

To further our discussion we turn to results from \cite{BF2} and
\cite{F1}.  Let $Y \rightarrow Z$ be a family of closed, spin 
Riemannian 2-manifolds, and $E \rightarrow Y$ a real virtual 
vector bundle with orthogonal connection $\nabla^E$.  Let $\Line{}{} = 
\Pfaff^{-1}_{Y/Z, E} \rightarrow Z$ be the inverse Pfaffian 
line bundle over the family of twisted Dirac operators.
Denote $I = [0,1]$ and take a path $\gamma : I \rightarrow Z$.  We 
form the pullback $\gamma^*Y \rightarrow I$ and observe that 
$\gamma^*Y$ is 3-manifold with boundary $\d (\gamma^*Y) = 
Y_{\gamma (1)} \sqcup -Y_{\gamma (0)}$.  If we place the standard 
metric on $I$, then we get a metric on $\gamma^*Y$ thanks to the metric 
$g_{\gamma^*Y /I}$ and the distribution of horizontal planes.  
From this data we get a twisted Dirac operator $D_{E}$ on $\gamma^*Y$ 
and so a $\rttau$-invariant
\begin{equation}\label{para trans}
\rttau_{\gamma^*Y}(D_E) \in
\Line{}{}_{\gamma (1)} \otimes
\Line{-1}{}_{\gamma (0)}.
\end{equation}
Let us look at the algorithm we have generated.  To a path in $Z$ -- 
using the $\rttau$-invariant -- we assign a linear map between the 
lines over the endpoints of the path.  This is exactly what occurs
in parallel transport.  In fact, the following theorem confirms that
the algorithm described above is just that.

\begin{theorem}[\cite{F1}, Theorem 3.1]
Let $E \rightarrow Y$ be a real rank zero virtual vector bundle
with orthogonal connection over the family of closed, spin 
Riemannian 2-manifolds $Y \rightarrow Z$.  
\begin{enumerate}
\item If $\gamma : I \rightarrow Z$ is a path then the linear map
 \begin{equation}\label{pf hol}
\rttau_{\gamma^*Y}(D_E) :
\Line{}{}_{\gamma (0)} \longrightarrow 
\Line{}{}_{\gamma (1)}
\end{equation}
is parallel transport over $\gamma$ with respect to the natural
connection on $\Line{}{}$.
\item Let $\Omega$ denote the curvature of the natural connection 
on $\Line{}{} \rightarrow Z$.  Then
\begin{equation}\label{pf curv}
\Omega = -\pi i \left[ \int_{Y/Z} 
\widehat{A}(\Omega^{Y/Z}) ch (\Omega^E) \right]_{(2)}.
\end{equation}
\end{enumerate}
\end{theorem}     

If $\gamma$ is closed path in $Z$ then we can
form a 3-manifold $Y_{\gamma}$ fibered
over the circle and $\rttau_{Y_{\gamma}}(D_E)$ is a 
unitary number.  Recall that $S^1$ has two spin structures: the 
{\it non-bounding} spin structure is the trivial double cover of the
circle and the {\it bounding} spin structure is the non-trivial 
double cover.  In this case we have the following theorem.

\begin{theorem}[\cite{BF2}, Theorem 3.18]
Suppose $\gamma : I \rightarrow Z$ is a closed path (which is
constant on $[0, \delta]$ and $[1- \delta,1]$ for some $\delta$).
Then the holonomy around $\gamma$ is given by
\begin{equation}\label{loop hol}
\text{\rm hol}_{\gamma} =
\begin{cases}
(-1)^{\text{\rm ind}_2 (D_{Y,E})}
\rttau_{Y_{\gamma}}(D_E), & 
\text{nonbounding spin structure on} S^1; \\
\rttau_{Y_{\gamma}}(D_E), & 
\text{bounding spin structure on} S^1.
\end{cases}
\end{equation}
\end{theorem}  

We point out some subtleties that have been sidestepped.  What 
we have written here is a really a specialized version of what appears
in \cite{F1} which applies more generally for $E$ of any rank
(and, as usual, in dimensions 2 (mod 8)).  By restricting $E$ to have
rank zero we avoid having to take an {\it adiabatic limit} of 
$\rttau$-invariants.  In general, upon taking the adiabitic limit,
the linear map \eqref{para trans} becomes independent of the
parametrization of the path $\gamma$.  The adiabatic limit is taken
over the metric we place on $I$; but on 3-manifolds, when $E$
is rank zero, we have seen that the $\rttau$-invariants are 
independent of the metric.  Thus our linear map is independent
of the parametrization of the path even before taking an adiabatic 
limit.  For details we defer to the source cited.

\subsection{The line bundle and symplectic structure}

Though we want our physical parameters to be the $G$-connections 
over a closed spin 2-manifold $Y$, to be completely egalitarian 
we allow the Riemannian metric to vary as well.
Just as the metrics are ``unphysical'' in the Lagrangian field
theory, so they are in the Hamiltonian field theory.  We explain this
shortly.  Thus, if we fix a principal $G$ bundle $Q \rightarrow Y$,
our {\it in nomine} parameter space is  $\Met{(Y)} \times 
\Cat{}{Q}$ but our {\it de facto} parameter
space is $\Cat{}{Q}$.  

Given a real virtual representation $\rho$,
we have a canonical family of twisted Dirac operators 
associated to $\Met{(Y)} \times \Cat{}{Q}$: 
$(g, A) \mapsto D_{g, \rho A}$; and so we have a canonical 
inverse Pfaffian line bundle 
\begin{equation}\label{canon line}
\Line{\rho}{(Q)} \rightarrow \Met{(Y)} \times \Cat{}{Q} .
\end{equation} 
This line bundle is the central figure in the Hamiltonian field theory 
over $Y$, just as $\rttau_X (D_{\rho})$ was in the Lagrangian field 
theory over the spin 3-manifold $X$.  We use the theorems above to 
formulate some of the geometric properties of $\Line{\rho}{(Q)}$. 

\begin{proposition}
Let $\rho$ be a real, rank zero virtual representation of a compact
Lie group $G$ and let $Q \rightarrow Y$ be principal $G$-bundle
over a closed spin 2-manifold.   
\begin{enumerate}
\item If $\gamma : I \rightarrow \Met{(Y)} \times \Cat{}{Q}$ 
is a path then 
\begin{equation}\label{hol}
\text{\rm hol}_{\gamma} = \rttau_{Y \times I}
 (D_{g^{\gamma}, \rho A^{\gamma}})
\end{equation} 
where $\text{\rm hol}_{\gamma}$ is the holonomy along $\gamma$ with
respect to the natural connection on $\Line{\rho}{(Q)}$ and 
$g^{\gamma}$, $ A^{\gamma}$ are the metric and $G$-connection on 
$Y \times I$ determined by $\gamma$.
\item The natural connection on $\Line{\rho}{(Q)}$ is flat 
along $\Met{(Y)}$ and if $\Omega$ denotes the curvature of 
$\Line{\rho}{(Q)}$ then 
\begin{equation}\label{curv}
\Omega (\dot{A_1}, \dot{A_2}) = 2 \pi i \int_Y 
\fpairing{\dot{A_1}}{\dot{A_2}}_{\rho}.
\end {equation}
where $\dot{A_j} \in \Omega^1(\text{\rm ad} Q)$, $j=1,2$ are tangent
vectors along $\Cat{}{Q}$.
\end{enumerate}
\end{proposition}

\begin{remark}
We explain in a more precise way what we mean when we 
say that the metric is ``unphysical'' in the Hamiltonian theory.
Since the connection is flat along the convex space 
$\Met (Y)$ we can use parallel transport to identify the 
fiber $\Line{\rho}{(g_1,B)}$ with the fiber 
$\Line{\rho}{(g_2,B)}$ in a way that is independent 
of the path in $\Met (Y) \times \{ B \}$ connecting the 
metrics $g_1$,$g_2$.  Also, since the curvature is flat along 
$\Met (Y)$ and does not depend on the metric in any way, $\Omega$ 
must lie in the image of the pullback by the projection $\Met{Y} 
\times \Cat{}{Q} \rightarrow \Cat{}{Q}$.  
In fact, assuming the pairing $\pairing{}{}_{\rho}$ is 
non-degenerate, $\omega = \Omega / 2 \pi i$ defines a symplectic form 
$\omega$ on $\Cat{}{Q}$.  From now on we require that the 
pairing be non-degenerate, so that our phase space is properly 
identified as the symplectic space $(\Cat{}{Q}, \omega)$. 
\end{remark}

\begin{proof}
Let $g_t$, $A_t$, $t = (t_1,t_2) \in [-1,1]^2$ be a two parameter
family of metrics and $G$-connections such that 
$(\d A_t /\d t_j)|_{t=0} = \dot{A_j}$.
This gives us a natural metric $g$ and connection $A$ on 
$Q \times [-1,1]^2$.  From \eqref{pf curv} we have
\begin{align}
\Omega = & -\pi i \left[ \int_Y 
\widehat{A}(\Omega^g) ch (\Omega^{\rho A}) 
\right]_{(2)} \nonumber \\
= & -\pi i \left[ \int_Y 
\fpairing{\Omega^{A}}{\Omega^{A}}_{\rho}
\right]_{(2)}. \label{curv1}
\end{align}
where the second equality follows from the unraveling of 
the characteristic polynomials $\widehat{A}$ and $ch$
and the fact that $\rho$ is rank zero.  If $\Omega^{A_t}$
denotes the curvature of the connection $A_t$ then,
$$
\Omega^A = dt_1 \wedge \dot{A_1} + 
dt_2 \wedge \dot{A_2} + \Omega^{A_t}.
$$
Plugging this into \eqref{curv1} we obtain the result.
\end{proof}

\subsection{The line bundle and symplectic reduction}

In the Lagrangian field theory over a spin 3-manifold $X$ we insisted 
that two $G$-connections be physically equivalent if they are isomorphic 
by some element of $\Mor{G}{X}$.  Indeed, we saw that the action 
$\rttau$ is invariant with respect to this symmetry.  The next few
propositions offer some consequences of this symmetry in the Hamiltonian
field theory.   

\begin{proposition}
The action of $\Mor{}{Q}$ on $\Cat{}{Q}$ lifts to 
$\Line{\rho}{(Q)}$ and the lifted action preserves the metric and 
the connection.  The induced moment map is 
\begin{equation} \label{gauge mom map}
\mu_{\zeta}(B) = \int_Y 
\fpairing{\Omega^B}{\zeta}_{\rho},
\end{equation}
where $\zeta \in \Omega^0(\text{\rm ad}Q)$ is an infinitesimal 
gauge transformation.
\end{proposition}

\begin{remark}
In the language of symplectic geometry, this proposition
implies that the moduli space of flat $G$-connections on which 
$\Mor{}{Q}$ acts freely is the symplectic quotient 
$\Cat{}{Q} // \Mor{}{Q}$.  As we shall explain in the proof,
there is an induced line bundle $\Line{\rho}{(Q)} \rightarrow
\Mod{}{Q}$ with a metric and a connection whose curvature is 
$2 \pi i$ times the symplectic form that naturally comes with the 
symplectic quotient.  
of flat connections 
\end{remark}

\begin{proof}
Recall that \eqref{isom on L} is a lift of the $\Mor{}{Q}$ action on 
$\Cat{}{Q}$ to $\Line{\rho}{(Q)}$ that preserves the metric.  
Given the formulation \eqref{hol} of parallel transport, \eqref{isom 
on rttau} implies that this lift preserves the connection.  Thus the action  
of $\Mor{}{Q}$ preserves the curvature of $\Line{\rho}{(Q)}$ and so 
the symplectic form $\omega$.  We can, therefore, compute the 
{\it moment map} of the symplectic action of $\Mor{}{Q}$.  

Before we compute the moment map of the $\Mor{}{Q}$ action
we recall how this is done in the context of automorphisms on line
bundles.  Whenever $\Line{}{} \rightarrow M$ is a hermitian 
line bundle with unitary connection and $\beta : G \rightarrow 
\text{Aut}(\Line{}{})$ is a $G$ action on $\Line{}{}$  that 
preserves the metric and connection, the moment map of $G$ 
action on $M$ is 
\begin{equation}\label{G mom map}
\mu_{\zeta}(m) = 
\frac{\text{vert}(\dot{\beta}(\zeta)_{u})}{2 \pi i},
\quad \zeta \in \Lie{g},
\end{equation}     
where $u \in \Line{}{}_m$ is a unitary element, $\dot{\beta}(\zeta)$
is the vector field on $\Line{}{}$ corresponding to $\zeta \in \Lie{g}$,
and $\text{vert}(\cdot)$ is the vertical part of the vector with respect
to the connection on $\Line{}{}$.  This is precisely the obstruction
to the connection descending to the quotient $\Line{}{}/G$.  If it
dissappears (as it does in our case for flat $G$-connections ) then 
the connection descends to the quotient bundle, as claimed in the
remarks following the proposition. 

Now we apply this to our situation.  Suppose $\zeta \in 
\Omega^0(\text{ad}Q)$ and $\phi_s \in \Mor{}{Q}$ is a path of gauge
transformations with $\phi_0 = id_Q$ and $\dot{\phi}_0 = \zeta$.  
Consider the path of $G$-connections $B_s = \phi_s^*B$ which forms 
 a connection $A_t$ on $Q \times [0,t]$.  To compute the vertical action we
``divide'' the  automorphsim $(\rho \phi_t)^*$ by the parallel transport 
$\rttau_{Y \times [0,t]}(D_{\rho A_t})$.  Then \eqref{glue rttau}
and \eqref{hol} imply that this number is the $\rttau$-invariant of the
$G$-connection $({\bf Q}_t, {\bf A}_t)$ over $Y \times S_b^1(t)$ gotten 
by gluing together the endpoints of $(Q \times [0,t], A_t)$ with $\phi_t$.  
Here $S^1_b(t) = [0,t]/ \{ 0,t \}$ denotes the circle of length $t$ with the
bounding spin structure.

To ease our computation, let $\boldsymbol{\phi}$ denote the automorphism of 
$Q \times [0,t]$ given by $(p,s) \mapsto (\phi_s^{-1}p, s)$.  It descends
to ${\bf Q}_t$ and we denote the descendent automorphsim by 
$\boldsymbol{\phi}$, as well.  A simple computation shows that
$$
\boldsymbol{\phi}^*{\bf A}_t = B - t \zeta \cdot \theta
$$     
where $\theta$ is the standard unit measure on $S^1$.
Since $\rttau$ is invariant with respect to gauge transformations,
we have reduced the computation of the vertical action to computing
$$
\rttau_{Y \times S_b^1}(D_{\rho (B - t \zeta \cdot \theta)}).
$$
An easy application of \eqref{d wrt conn} gives us the infinitesimal
vertical action
$$
\text{vert}(\dot{(\rho \phi_t^*)}|_{t=0}) = 
\int_Y \fpairing{\Omega^B}{\zeta}_{\rho}
$$
and this, together with \eqref{G mom map}, gives us 
\eqref{gauge mom map}. 
\end{proof}

The gluing law \eqref{glue rttau} (perhaps more appropriately 
called a ``cutting law'') is used to good effect in the proof above.
We use it again to obtain a result for the action of the stabilizer 
subgroup $\Mor{}{B} \subset \Mor{}{Q}$.  In general, the 
gauge transformations do not act freely on $\Cat{}{Q}$.  The
subgroup $\Mor{}{B}$ are the elements of 
$C^{\infty}(\text{Ad}Q)$ that are parallel with respect 
to $B$.  Therefore, the value of a stabilizer of $B$ at $y \in Y$ 
-- which lies in $\Aut (Q_y)$ -- commutes with the holonomy group
of $B$ at $y$.  

\begin{proposition}\label{stabilizer}  
The action of $\Mor{}{B}$ on $\Line{\rho}{(B)}$ is constant
on components of $\rho \Mor{}{B}$ and so factors through an
action of the finite group $\pi_0 (\Mor{}{B})$ on 
$\Line{\rho}{(B)}$. 
\end{proposition}
 
 \begin{proof}
 Let $\phi_t$, $t \in I$ be a path in  $\Mor{}{B}$ and let
 $({\bf Q}_t, {\bf B}_t )$ denote the vector bundle with
 connection over on 
 $S_b^1 \times Y$ gotten by gluing the together the endpoints
 of $(Q \times I, B )$ with $\phi_t$.  In this way we
 obtain a connection over $I \times S_b^1 \times Y$
 which is flat along $I \times S^1$.  We use parallel transport
 along $I$ to obtain a morphism between the vector bundles with 
 connection $({\bf Q}_0, {\bf B}_0)$ and  
 $({\bf Q}_1, {\bf B}_1)$.  Thus 
 $$
 \rttau_{S_b^1 \times Y}(D_{\rho {\bf B}_0 }) =
 \rttau_{S_b^1 \times Y}(D_{\rho {\bf B}_1 })
 $$
 and so, according to \eqref{glue rttau}, $\phi_0$ and $\phi_1$ 
 induce the same automorphism on $\Line{\rho}{B}$.
 \end{proof}
 
\subsection{The functor of the line bundle}

We end this chapter with an observation of how the 
Lagrangian and Hamiltonian field theories are related to one
another.  Let $X$ be any compact spin 3-manifold with boundary 
$\d X$.  There is a restriction functor 
$r_X : \Cat{G}{X} \rightarrow \Cat{G}{\d X}$ which 
sends $(P,A) \mapsto (P,A)|_{\d X}$.  This functor restricts
to flat connection and descends to equivalence classes so 
that we have the commutative diagram
$$\begin{CD}
\Mod{G}{X}  @>>>  \overline{\Cat{G}{X}} \\
@V{r_X}VV                               @VV{r_X}V \\
\Mod{G}{\d X}  @>>>   \overline{\Cat{G}{X}}.
\end{CD}$$
Having said that we offer the following proposition.

\begin{proposition}
The map $r_X : \Mod{G}{X} \rightarrow \Mod{G}{\d X}$
is Lagrangian.  In fact, the action $\rttau_X$ is a 
flat section of the pullback $r_X^* \Line{\rho}{(\d X)} \rightarrow
\Mod{G}{X}$ so that the induced symplectic form $r_X^* \omega$
vanishes.
\end{proposition}
  
This proposition implies that the holonomy on $r_X^* \Line{\rho}{(\d X)}$
is trivial.  The image of $r_X$ is an example of a {\it Bohr-Sommerfeld 
orbit} of the line bundle $\Line{\rho}{(\d X)}$.  These play a pivotal role
in quantization with {\it real polarizations}.  We will have more to 
say about this when we quantize the Hamiltonian theory over a genus one 
surface.

\begin{proof}
To simplify the exposition we restrict to the connections on a 
particular $G$-bundle $P \rightarrow X$.
Let $\Cat{}{P}^{\text{flat}}$ denote the subset of flat connections
on $P$.  Over this space the variation formula \eqref{d wrt conn2}
gives
$$
\nabla_{\dot{A}} \rttau_X = \pi i \cdot
\rttau_X \int_X \fpairing{\dot{A}}{\Omega^A}_{\rho}
$$ 
where $\dot{A}$ is some vector at $T_A \Cat{}{Q}$.  Since
$A$ is flat the right hand side is zero and so $\rttau_X$
is a flat section.  

It remains to show that, if $A$ and $A|_{\d X}$ represent smooth 
points of $\Mod{G}{X}$ and $\Mod{G}{\d X}$ respectively, then
$$
2 \dim \text{ image}( (r_X)_* ) = \dim H^1(X; d_A).
$$
As we have nothing new to add we defer to the proof
that appears in Proposition 3.27 of \cite{F2}.
\end{proof}

We end by summarizing the main point of the Hamiltonian field
theory, that being the assignment
$$
Y \mapsto ( \Line{\rho}{(Y)} \rightarrow \Mod{G}{Y} ).
$$
This clearly obeys the same functoriality, orientation, additivity and 
gluing laws that appear in Proposition \ref{functorial2}.  The 
Hamiltonian theory makes contact with the Lagrangian theory
when $Y = \d X$.  In this case we have the assignment
$$
X \mapsto (\rttau_X : \Mod{G}{X} \rightarrow 
r_X^*\Line{\rho}{(\d X)}), 
$$ 
and we have seen this adheres to the same four rules mentioned
above.  This ends our study of the classical spin-Chern-Simons
field theory. 

\section{Appendix}

\subsection{$KO$ of a compact 3-manifold}

We compute the $KO$ group of a closed, connected 3-manifold 
in terms of its $\zmod2 $-valued cohomology.  First of all we 
have that, for $X$ connected, 
$$
KO(X) \cong \Z \oplus \widetilde{KO}(X)
$$ 
where $\widetilde{KO}(X)$ is the kernel of the rank map.

Now we consider the group structure on 
$H^1(X;\zmod2 ) \times H^2(X;\zmod2 )$ given by the
product
$$
(a_1, b_1) \cdot (a_2, b_2) = (a_1 + a_2, a_1 \smile a_2 + b_1 + b_2).
$$
We let $H^1(X;\zmod2 ) \ltimes H^2(X;\zmod2 )$ denote the set 
$H^1(X;\zmod2 ) \times H^2(X;\zmod2 )$ with this group structure. 

\begin{theorem}\label{ko of x}
For $X$ a closed, connected, compact 3-manifold, the map 
\begin{align*}
\widetilde{KO}(X) \longrightarrow &
H^1(X;\zmod2 ) \ltimes H^2(X;\zmod2 ) \\ 
E - F \longmapsto &
(w_1(E) + w_1(F) ,
 w_2(E) + w_2(F) +  w_1(E) w_1(F) + w_1(F)^2)  
\end{align*}
is an isomorphism. 
\end{theorem}

To prove Theorem \ref{ko of x} we will use the following fact 
about oriented vector bundles over a compact 3-manifold.

\begin{proposition}\label{w2 classifies}
For $n \geq 3$ the 2nd Stiefel-Whitney class provides a 
one-to-one correspondence between topological $SO_n$ bundles
over $X$ and  classes in $H^2(X ; \zmod2 )$. 
\end{proposition} 

\begin{proof}
We first show injectivity.  From the short exact sequence 
$$
1 \rightarrow \zmod2 \longrightarrow 
Spin_n \longrightarrow SO_n \rightarrow 1
$$
we get an exact sequence of \v{C}ech cohomology groups
$$
\dots \rightarrow H^1(X;Spin_n) \longrightarrow
H^1(X;SO_n) \xrightarrow{w_2} 
H^2(X ; \zmod2 )
$$
For $n \geq 3$ $Spin_n$ is simply connected so that, according to
a simple argument in obstruction theory, any $Spin_n$ principal 
bundle over a 3-manifold is topologically trivial.  Thus the map
$w_2$ is injective. 

We now show $w_2$ is surjective.  To see this one takes the short exact 
sequence
$$
0 \rightarrow \Z \xrightarrow{2 \times} 
\Z \xrightarrow{\text{mod 2}} 
\zmod2 \rightarrow 0
$$
and the induced long exact sequence in cohomology
$$
\dots \rightarrow H^2(X; \Z) 
\xrightarrow{\text{mod 2}} H^2(X; \zmod2) 
\xrightarrow{\beta}  
H^3(X; \Z) \rightarrow \dots .
$$
The image of $\beta$ in $H^3(X; \Z) \cong \Z$ will be torsion and 
so must be zero.  Thus the ``mod 2'' cohomology map is surjective.  
Of course, $H^2(X ; \Z)$ parametrizes the topological $SO_2$ bundles
over $X$ via the 1st Chern class 
$$
c_1: H^1(X; SO_2) \longrightarrow H^2(X; \Z ).
$$  
The standard inclusion homomorphism $SO_2 \hookrightarrow SO_n$ 
and the associated bundle map 
\begin{align*}
H^1(X ; SO_2) & \longrightarrow H^1(X ; SO_n) \\
[P]                   & \longmapsto  [P \times_{SO_2} SO_n]
\end{align*}
give us the commutative diagram
$$
\begin{CD}
H^1(X ; SO_2)   @>\text{assoc. bundle}>>  H^1(X;SO_n)  \\
@V{c_1\text{ mod 2}}VV                                         @VVw_2V \\
H^2(X; \zmod2)   @>\text{identity}>>              H^2(X; \zmod2)
\end{CD} 
$$ 
so that, since the ``mod 2'' map is surjective, $w_2$ is surjective. 
\end{proof}

With this result at our disposal we can now prove Theorem \ref{ko of x}.

\begin{proof}[Proof of Theorem \ref{ko of x}]
That the map is well-defined and a homomorphism easily 
follow from the properties of the Stiefel-Whitney classes.  

To prove surjectivity, choose any $b \in H^2(X ; \zmod2)$ and 
$a \in H^1(X ; \zmod2)$.  Let $\ell$ be a real line bundle
over X such that $w_1(\ell ) = a$; and, according to Lemma
(A.1.2), we can choose an oriented rank 2 bundle $E$ such that
$w_2(E) = b$.  Then 
$$
(w_1,w_2)(E \oplus \ell) = (a,b)
$$ 
proving surjectivity.

To prove injectivity we consider the kernel of the homomorphism.
Indeed, suppose $E - F \in \widetilde{KO}(X)$ is such that
\begin{align*}
w_1(E)  & = w_1(F) \quad \text{and} \\
w_2(E)  & = w_2(F) + w_1(F)^2 + w_1(F) w_1(E)
\end{align*}
The first equality implies $\Det E \cong \Det F$ so that in $KO(X)$
$$
E - F = (E \oplus \Det E) - (F \oplus \Det F) 
$$
Since the map is well defined and
$$
w_1(E \oplus \Det E) = w_1(F \oplus \Det F) = 0
$$
we must have that
$$
w_2(E \oplus \Det E) = w_2(F \oplus \Det F)
$$
From Lemma \ref{w2 classifies} we know this equality implies the
equivalence 
$$
E \oplus \Det E \cong F \oplus \Det F
$$ 
so that in $KO(X)$
$$
E - F = (E \oplus \Det E) - (F \oplus \Det F) = 0.
$$
Thus the kernel of the homomorphism is trivial.  
\end{proof}

\subsection{Cobordism groups}

In this section we compute certain cobordism groups that 
pop up Chapter 2.  In particular, they appear in the proofs
to Propositions \ref{properties of q} and \ref{level}.

\begin{proposition}\label{cobordism grps}
Let $\Omega^{\text{\rm spin}}_n (M)$ denote the degree-$n$ 
spin cobordism group of the topological space $M$.  Then
$$
\Omega^{\text{\rm spin}}_3 (BSO \wedge B(\zmod2)) = \zmod2
\quad \text{\rm and} \quad
\Omega^{\text{\rm spin}}_3 (BSO) \cong \zmod2.
$$
and the generator of $\Omega^{\text{\rm spin}}_3 (BSO)$ is 
represented by $S^1_{nb} \times F \rightarrow S^1_{nb} \times S^2$.
Here $F \rightarrow S^2$ is an oriented rank zero virtual vector 
bundle with non-trivial $w_2$ and $S^1_{nb}$ is the circle with the 
non-bounding spin structure.
\end{proposition}

\begin{proof}
The first isomorphism can be seen as follows.  Since $BSO$ is 
1-connected and $B(\zmod2)$ is connected the wedge product 
$BSO\wedge B(\zmod2)$ is 2-connected.  The Hurwitz theorem implies 
that the first non-trivial homology (for any coefficient group) is $H_3$ 
and from the $E^2$ term of the Atiyah-Hirzebruch spectral sequence it 
is clear we need only determine $H_3 (BSO\wedge B(\zmod2))$.  This 
is easily computed by considering the long exact homology sequence 
induced by the maps
$$
BSO \vee B(\zmod2) \hookrightarrow
BSO \times B(\zmod2) \longrightarrow
BSO \wedge B(\zmod2). 
$$ 
One then sees that $H_3 (BSO\wedge B(\zmod2)) \cong \zmod2$ and 
this proves the first isomorphism of the propositon.

The first isomorphism is used to prove Proposition \ref{properties of q}
and now we use that proposition to prove the second isomorphism.
From the $E^2$ term of the Atiyah-Hirzebruch spectral sequence we can
conclude that $\Omega^{\text{\rm spin}}_3 (BSO_N)$ is isomorphic to
either $\zmod2$ or $\{ 0 \}$ for any $N > 2$. 
We will now show that certain structures on $S^1 \times S^2$ cannot 
be the boundary of structures on a 4-manifold so that the cobordism
group must be isomorphic to $\zmod2$.

Let $(Q, B) \rightarrow S^2$ be an $SO_N$-connection such that 
$Q$ has non-trivial $w_2$. Contrary to the claim of the proposition, 
we assume that there exists $(P,A) \rightarrow M$ such that 
\begin{equation}\label{cont boundary}
\d ( (P,A) \rightarrow M ) = 
( S^1_{nb} \times (Q,B) \rightarrow S^1_{nb} \times S^2 ).
\end{equation}     
On the other hand, there is no question that there exists 
$(P',A') \rightarrow D^2 \times S^2$ such that
$$
\d ( (P',A') \rightarrow D^2 \times S^2 ) = 
( S^1_b \times (Q,B) \rightarrow S^1_b \times S^2 ).
$$
Let $\rho = id_{SO_N} - N$ be the identity representation minus
the rank $N$ trivial representation and let $\ell$ represent the non-trivial
element of $H^1(S^1;\zmod2)$.  According to Proposition \ref{properties of q}
$$
\frac { \rttau_{S^1_{nb} \times S^2}(D_{\rho B}) }
{ \rttau_{S^1_b \times S^2}(D_{\rho B}) } =
\rttau_{S^1_b \times S^2}(D_{\rho B \otimes (\ell - 1)}) =
(-1)^{w_2(Q) \smile \ell} = -1.
$$
However \eqref{cont boundary} implies
$$
\frac { \rttau_{S^1_{nb} \times S^2}(D_{\rho B}) }
{ \rttau_{S^1_b \times S^2}(D_{\rho B}) } =
\exp \pi i \int_{M \sqcup -(D^2 \times S^2)}
\fpairing{\Omega^{A \sqcup A'}}{\Omega^{A \sqcup A'}}_{\rho} =
1
$$
where the last equality follows from the fact the integral 
is an even integer (as it is the index of the twisted chiral
Dirac operator $D_{\rho (A \sqcup A')}$ on $M \sqcup
-(D^2 \times S^2)$).  Thus we get a contradiction and thereby
prove the second isomorphism and the claim following it.
\end{proof}

\subsection{Computations on the 3-Torus}

Here we compute the quadratic form $q$ as defined in 
in \eqref{quad form} when the closed 3-manifold is 
a 3-torus.  These particular computations are used in 
the cobordism argument for parts (2) and (3) of 
Proposition \ref{properties of q}.

\begin{proposition}\label{linearity}
Let $E \in \widetilde{KO}(X)$ be such that $w_1(E) = 0$.
Then $q_{\sigma }(E, \ell )$ is independent of the spin
structure $\sigma$ and depends linearly on $\ell \in H^1(X;\zmod2)$.
\end{proposition}

\begin{proof}
An easy computation shows that if $\ell_0 \in H^1(X;\zmod2)$ then 
$w_1(E \otimes \ell_0) = 0$ and $w_2(E \otimes \ell_0) = w_2(E)$.
Thus, by \ref{ko of x} $E \otimes \ell_0 = E$ in $KO(X)$.
This implies
$$
q_{\sigma + \ell_0}(E, \ell) = 
q_{\sigma }(E \otimes  \ell_0, \ell) =
q_{\sigma }(E, \ell)
$$
where the first equality follows from the fact that the Dirac operator
for $\sigma + \ell_0$ is just the Dirac operator for $\sigma $ twisted
by $\ell_0$.  This proves independence with respect to the spin 
structure.

Now consider $q_{\sigma }(E, \ell_1 \otimes \ell_2)$.  Using the
argument in the above paragraph we see that   
$$
q_{\sigma }(E, \ell_1 \otimes \ell_2) =
q_{\sigma }(E \otimes \ell_1, \ell_1 \otimes \ell_2) =
q_{\sigma }(E, \ell_1) \cdot q_{\sigma }(E, \ell_2)
$$
where the second equality follows easily from the definition
of $q$ \ref{quad form} and the fact that (in the case $w_1(E)=0$)
$q$ is $\zmod2^{\times }$-valued. 
\end{proof}

We are now in a good position to make some actual
computations.  In particlar we look at $q$ on the 
3-torus $T^3$.  What we find is that 

\begin{proposition}\label{q on 3torus}
Let $E \in \widetilde{KO}(T^3)$ be such that $w_1(E) = 0$.
Then for any $\ell \in H^1(T^3;\zmod2)$ 
$$
q_{\sigma }(E, \ell ) = (-1)^{w_2(E) \smile \ell}.
$$
\end{proposition}

\begin{proof}
Take $T^3 = \R^3 / \Z^3$ with the three standard projections
$\pi_j: T^3 \rightarrow T$ given by $(x_1, x_2, x_3) \mapsto x_j$
for $j=1,2,3$.  Let $\ell_0 \in H^1(T;\zmod2)$ be the non-trivial 
element.  Then the elements $\ell_j = \pi_j^* \ell_0$ generate 
$H^1(T^3;\zmod2)$ and the elements $\ell_j \smile \ell_i$
generate $H^2(T^3;\zmod2)$.  We prove the proposition for a 
set of generators and then appeal to linearity for the other cases.

Without loss of generality we assume $w_2(E) = \ell_1 \smile 
\ell_2$ and consider the cases $\ell = \ell_1$ or $\ell_2$.  
Then $w_2(E) \smile \ell = 0$.  On the other hand, all of the 
topological data is trivial along the 1-cycle $c_3 = (0 \oplus 0 \oplus \R )
/ \Z^3$.  We are free to choose the smooth parameters to be 
trivial along $c_3$ as well.  In fact, according to Proposition \ref{linearity}
just above, we can even assume that the spin structure $\sigma$ extends
across a disk bounded by $c_3$.
This allows us to implement the gluing formula \ref{glue rttau} on 
each of the $\rttau$ factors appearing in $q$.  We find
$$
\rttau_{T^2 \times S^1_b} (D_{E \otimes \ell}) = 1 
\quad \text{and} \quad
\rttau_{T^2 \times S^1_b} (D_{E}) = 1 
$$
so that, in this case,
$q_{\sigma }(E, \ell ) = 1 = (-1)^{w_2(E) \smile \ell}$.

Now we consider the case $\ell = \ell_3$.  Then $w_2(E) \smile \ell 
= 1$.  On the other hand, we can take all of the parameters to be 
trivial along $c_3$ (including a bounding spin structure) except for the
line bundle $\ell_3$.  But we can still use the gluing formala 
\ref{glue rttau}.  Indeed we see that 
\begin{align*}
\rttau_{T^2 \times S^1_b} (D_{E \otimes \ell_3}) = &
\rttau_{T^2 \times S^1_{nb}} (D_{E}) = (-1)^{ind_2 (D_{T^2, E'})} \\
\text{and} \quad
& \rttau_{T^2 \times S^1_b} (D_{E}) = 1 
\end{align*}
where $D_{T^2, E'}$ is the Dirac operator on $T^2$ twisted by the 
virtual vector bundle $E'$ which is the restriction of $E$ to $T^2 \subset
T^3$.  Since $ind_2 (D_{T^2, E'}) = w_2(E') = 1$ we get that, in this 
case, $q_{\sigma }(E, \ell ) = -1 = (-1)^{w_2(E) \smile \ell}$
\end{proof}

\subsection{The ``Level'' of the Theory}

To define the action of our spin-Chern-Simons theory we require an element 
of the virtual oriented representation ring, $RSO(G)$.
In fact, to eliminate metric dependence, we require that the element have 
rank zero.  Such elements form an ideal which we denote by
$\widetilde{RSO}(G)$.

A represention, $\rho: G \rightarrow SO_N$, generates a unique homotopy
class of maps, $B_{\rho} : BG \rightarrow BSO$, between classifying spaces.
Homotopy classes of maps into $BSO$ is a bit too rich in structure for our
needs.  In order to obtain a leaner structure we cap-off all of the homotopy
groups of $BSO$ above $\pi_4$.  What is left is a space whose only nontrivial
homotopy lies in $\pi_2$ and $\pi_4$.  In fact, this capped-off $BSO$ is
homotopic to a fibration, the total space of which we'll denote by $E^4$.  The
fibration is

$$ \begin{CD}
K(\Z,4)   @>\text{inclusion}>>   E^4 \\
@.                 @VVw_2V \\
            @.     K(\Z/2\Z,2)
\end{CD} $$

The map $w_2$ is an extension of the usual second Steiffel-Whitney class on
$BSO$.
Thus $E^4$ is a twisted product of Eilenberg-MacLane spaces.  Such fibrations
with
fiber $K(A,m)$ and base $K(B,n)$ are classified by elements of
$H^{m+1}(K(B,n),A)$.
The fibration seen above is determined by $\beta \circ Sq^2(\iota)$, where
$\beta$ is the
Bockstein homomorphism and $\iota$ is the fundametal class of
$H^2(K(\Z/2,2);\Z/2)$.
The ``levels'' in this theory are homotopy classes of maps from $BG$ into
$E^4$.  This has a group structure since $E^4$ is, much like an
Eilenberg-Maclane space, homotopic to a loop space.  Indeed, if we take the
fibration
$$ 
\begin{CD}
K(\Z,5)   @>\text{inclusion}>>   E^5 \\
@.                 @VVwV \\
            @.     K(\Z/2\Z,3)
\end{CD} 
$$
which is determined by $\beta \circ Sq^2(\iota')$ -- where $\iota'$ is the
fundamental
class of $H^3(K(\Z/2,3);\Z/2)$ -- we have that $E^4 \sim \Omega E^5$.

In fact, the two previous fibrations are part of something in homotopy theory
called a ``spectrum''.  For any counting number, $N$, we have a fibration
$$ 
\begin{CD}
K(\Z,N)   @>\text{inclusion}>>   E^N \\
@.                 @VVwV \\
            @.     K(\Z/2\Z,N-2)
\end{CD} 
$$
such that  $E^{N-1} = \Omega  E^N$.  Each such fibration is determined by 
$\beta \circ Sq^2(\iota_{N-2})$, where $\iota_{N-2}$ is the fundamental 
class of $H^{N-2}(K(\Z/2,N-2);\Z/2)$.   Such a spectrum of fibrations induces 
a long exact sequence of homotopy groups.  In particular we have
$$
\dots \rightarrow [X, K(\Z, N) ] \rightarrow [ X,E^N ] \rightarrow
[ X,K(\Z/2\Z, N-2) ] \rightarrow [ X, K(\Z,N+1)] \rightarrow \dots
$$
or more succinctly
$$
\dots \rightarrow H^N(X;\Z) \rightarrow E^N(X) \rightarrow
H^{N-2}(X;\Z/2\Z) \xrightarrow{\beta \circ Sq^2} H^{N+1}(X;\Z) 
\rightarrow \dots
$$
We will be interested in the following portion of the long exact sequence for
$BG$:
$$
\dots \rightarrow H^4(BG;\Z) \rightarrow E^4(BG) \rightarrow 
H^2(BG;\Z/2) \rightarrow \dots
$$
An easy argument shows that the first map is injective.  Indeed, the
map $Sq^2$ is zero on $H^1(BG;\zmod2)$ so that $\beta \circ Sq^2 : 
H^1(BG; \zmod2) \rightarrow H^4(BG;\Z)$ is zero as well.  From this
definition it is clear that when $G$ is simply connected $E^4(BG) \cong
H^4(BG;\Z)$ since $H^2(BG;\zmod2) = 0$.
 
The first Pontryagin map $p_1:BSO \longrightarrow K(\Z,4)$ extends to a map
on
$E^4$.  One can show that 

\begin{proposition}\label{p1}
The sequence of homomorphisms
$$ 
H^4(BG;\Z) \longrightarrow E^4(BG) 
\xrightarrow{p_1}  H^4(BG;\Z)
$$
is multiplication by $2$ in $H^4(BG;\Z)$.
\end{proposition}

\begin{proof}
We begin by noting that the fibration 
$$
K(\Z , 4) \hookrightarrow E^4 
\xrightarrow{w_2} K(\zmod2 , 2) 
$$
can be derived from the fibration
$$
BSpin \hookrightarrow BSO 
\xrightarrow{w_2} K(\zmod2 , 2) 
$$
by capping off all the homotopy generators of $\pi_n(BSpin)$ for
$n > 4$.
Thus, when we consider the sequence of maps 
$$
K(\Z ,4) \hookrightarrow E^4 
\xrightarrow{p_1} K(\Z ,4)
$$
we are effectively comparing the pullback of $p_1 \in H^4(BSO;\Z)$ to a
generator of $H^4(BSpin;\Z)$.  If we denote the pullback by $p_1$ as well
(as is commonly done in the literature) then it is well known that there
is a generating class ``$p_1/2$'' $\in H^4(BSpin;\Z)$ such that $p_1 = 2 
\cdot (p_1/2)$.  This proves the proposition.  

\end{proof}

\begin{proposition}\label{no torsion}
If $G$ is a connected compact lie group then $H^4(BG;\Z)$ is torsionless.
\end{proposition}
 
\begin{proof}
Let $T \subset G$ be a maximal torus.  Consider the cohomology 
spectral sequence $E_r$ of the fibration 
$G/T \hookrightarrow BT \rightarrow BG$.
We have two useful facts to help us along: 
$H^{\text{odd}} (G/T) = 0$ and 
$H^{\text{even}} (G/T)$ is torsionless.

Now $E_2^{4,0} = H^4(BG;\Z)$ and its clear from all of the zeros 
in the degree 3 diagonal line that $E_{\infty }^{4,0} = E_2^{4,0}$.
Since $E_{\infty }^{4,0}$ is isomorphic to a subgroup of $H^4(BT)$,
which is torsionless, we conclude that $H^4(BG;\Z)$ is torsionless.
\end{proof} 

\begin{proposition}\label{injective}
The homomorphism
$$
(p_1,w_2):E^4(BG) \rightarrow H^4(BG;\Z) \oplus H^2(BG;\Z/2)
$$
is injective up to 2-torsion elements of $H^4(BG;\Z)$ so that, if $G$ is 
connected, $(p_1,w_2)$ is injective.
\end{proposition}

As an easy corollary we have

\begin{corollary}\label{group struct}
Let $\lambda(\rho)$ be the class in $E^4(BG)$ induced by the represention
$\rho$.
Then we have $\lambda(\rho_0 \oplus \rho_1) = \lambda(\rho_0) +
\lambda(\rho_1)$.
\end{corollary}

As an example -- and because we will later use these results when we consider the
quantum theory -- we consider $E^4(BSU_2)$ and $E^4(BSO_3)$.  Since $SU_2$ 
is simply connected we easily have 
$$
E^4(BSU_2) \cong H^4(BSU_2;\Z) = \Z \cdot c_2
$$
where $c_2$ is the 2nd Chern class.  More relevant to spin-Chern-Simons
is the image of the map 
$$
\lambda : \widetilde{RO}(BSU_2) 
\longrightarrow E^4(BSU_2).
$$ 
Let $\rho : SU_2 \rightarrow SO_4$ be the realization of 
the standard $SU_2$ representation on $\C^2$.  We take the rank zero
representation $(\rho - 4)$; that is, $\rho$ minus the trivial 4 dimensional
representation.  With some foresight, we denote ${\bf 1'} 
= \lambda (\rho - 4)$ and claim that $E^4(SU_2) = \Z \cdot {\bf 1'}$. 
Indeed, $p_1 ({\bf 1'}) = -2c_2$; and combined with Proposition \ref{p1} 
this proves the claim.  Then, according to Corollary 
\ref{group struct}, to hit all of the levels of $SU_2$ spin-Chern-Simons we 
need only consider integer multiples of $(\rho - 4)$.

Of course, $SO_3$ is not simply connected.  In fact, $H^2(BSO_3;\zmod2) \cong 
\zmod2$ is generated by the 2nd Stiefel-Whitney class.  We 
consider the image of the map 
$$
\lambda : \widetilde{RO}(BSO_3) 
\longrightarrow E^4(BSO_3).
$$
In particular we consider the element ${\bf 1} = \lambda( id_{SO_3} - 3)$ 
and claim that $E^4(BSO_3) = \Z \cdot {\bf 1}$.  Indeed,  $w_2 ({\bf 1}) \in
H^2(BSO_3;\zmod2)$ is the non-trivial element so that $w_2$ is surjective. 
Thus we have a short exact sequence
$$
0 \rightarrow (H^4(BSO_3;\Z) \cong \Z ) \rightarrow E^4(BSO_3)
\xrightarrow{w_2} ( H^2(BSO_3;\zmod2) \cong \zmod2 ) \rightarrow 0
$$
where $H^4(BSO_3;\Z)$ is generated by the 1st Pontryagin class.  This implies
that $E^4(BSO_3)$ is isomorphic is either $\Z \oplus \zmod2$ or $\Z$.  If
it is isomorphic to the former then, with respect to that isomorphism, 
$\text{image} (p_1 \oplus w_2) = 2\Z \oplus \zmod2$.  However, it is clear that
$$
(p_1 \oplus w_2)({\bf 1}) = 
1 \oplus (1 \quad (\text{mod} 2) ),
$$ 
and so it must be that $E^4(BSO_3)$ is isomorphic to $\Z$.  
Now Propositions \ref{p1} and \ref{group struct} imply that ${\bf 1}$ is a 
generator.

We end these considerations by pointing out that the standard 2:1 covering
homomorphism $\beta : SU_2 \rightarrow SO_3$ induces the homomorhpism
\begin{align*}
B_{\beta}^* : E^4(BSO_3) & \longrightarrow E^4(BSU_2) \\
                       k \cdot {\bf 1} & \longmapsto  2k \cdot {\bf 1'}.
\end{align*}
Indeed, consider the commutative diagram 
\begin{equation}\label{so3 su2 comm diag}
\begin{CD}
H^4(BSO_3;\Z)   @>>>   E^4(BSO_3) @>w_2>> H^2(BSO_3;\zmod2) \\
@V{B_{\beta}^*}VV    @V{B_{\beta}^*}VV   @VV{B_{\beta}^*}V  \\
H^4(BSU_2;\Z)   @>>>   E^4(BSU_2) @>w_2>> 0 .
\end{CD}
\end{equation}
A simple argument in Chern-Weil theory shows that $B_{\beta}^*$ sends the 
generator of $H^4(BSO_3;\Z)$ to 4 times the generator of $H^4(BSU_2;\Z)$.  Thus,
the diagram \eqref{so3 su2 comm diag} is isomorphic to the diagram 
$$
\begin{CD}
\Z   @>{\times 2}>>   \Z @>{(\text{mod} 2)}>> \zmod2 \\
@V{\times 4}VV    @V{B_{\beta}^*}VV   @VV0V  \\
\Z   @>{\times 1}>>   \Z @>0>> 0 
\end{CD}
$$
and from this diagram's commutativity it is clear that $B_{\beta}^*$ is 
effectively multiplication by 2.

\bibliographystyle{plain}

\end{document}